\documentclass[11pt]{amsart}

\allowbreak
\usepackage{amsfonts,amssymb,amscd,amsmath,latexsym,amsbsy}
\usepackage{graphicx}
\usepackage{amssymb}
\usepackage{amsfonts}
\usepackage{latexsym}
\usepackage{mathrsfs}
\usepackage{mathtools}
\usepackage{amsthm}
\usepackage{listings}
\usepackage{comment}
\usepackage{tikz-cd}
\usepackage{titlesec}
\usepackage{leftidx}

\allowdisplaybreaks
\newtheorem{thm}{Theorem}[section]
\newtheorem{defn}[thm]{Definition}
\newtheorem{lmn}[thm]{Lemma}
\newtheorem{prop}[thm]{Proposition} 
\newtheorem*{rmk*}{Remark} 
\newtheorem{rmk}[thm]{Remark}

\usepackage{mathtools}

\numberwithin{equation}{section}

\newcommand{\R}{\ensuremath{\mathbb{R}}}

\newcommand{\F}{\ensuremath{\mathsf{F}}}

\newcommand{\supp}{\ensuremath{\mathrm{supp}}}

\newcommand{\la}{\ensuremath{\langle}}
\newcommand{\ra}{\ensuremath{\rangle}}

\newcommand{\Id}{\ensuremath{\mathrm{Id}}}

\titleformat{name=\section}{}{\thetitle.}{0.8em}{\centering\scshape}
\titleformat{name=\subsection}[runin]{}{\thetitle.}{0.5em}{\bfseries}[.\\]
\titleformat{name=\subsubsection}[runin]{}{\thetitle.}{0.5em}{\itshape}[.\\]
\titleformat{name=\paragraph,numberless}[runin]{}{}{0em}{}[.]
\titlespacing{\paragraph}{0em}{0em}{0.5em}
\titleformat{name=\subparagraph,numberless}[runin]{}{}{0em}{}[.]
\titlespacing{\subparagraph}{0em}{0em}{0.5em}

\title{The tensorial X-ray transform on asymptotically conic spaces}
\author{Qiuye Jia}
\author[Andras Vasy]{Andr\'as Vasy}
\thanks{The authors gratefully acknowledge support from the National
  Science Foundation under grant number DMS-1953987.}

\date{\today}

\address{Department of Mathematics, Stanford University, Stanford, CA
94305-2125, U.S.A.}
\email{andras@math.stanford.edu}
\email{jqy@stanford.edu}

\begin{document}

\begin{abstract}
In this paper we show the invertibility of the geodesic X-ray
transform on one forms and 2-tensors on asymptotically conic
manifolds, up to the natural obstruction, allowing existence of
certain kinds of conjugate points. We use the 1-cusp
pseudodifferential operator algebra and its semiclassical foliation
version introduced and used by Vasy and Zachos, who showed the same type
invertibility on functions.

The complication of the invertibility of the tensorial X-ray
transform, compared with X-ray transform on functions, is caused by
the natural kernel of the transform consisting of `potential
tensors'. We overcome this by arranging a modified solenoidal gauge condition, under which we have the invertibility of the X-ray transform.
\end{abstract}

\maketitle
\section{Introduction}
For $n \geq 2$, the geodesic X-ray transform $I$ on a $n-$dimensional Riemannian manifold $(M,g)$, possibly with boundary, of a rank $m$ tensor $f$ is defined by
\begin{align}
If(\gamma) = \int \la f(\gamma(s)),\dot{\gamma}^m(s) \ra ds , \label{xray_defn}
\end{align}
where the paring is given by $\la f,v^m \ra =
\sum_{i_1,i_2,...,i_m}f_{i_1i_2...i_m}v^{i_1}v^{i_2}...v^{i_m}$ in
local coordinates (or a local frame), and we assume conditions on $\gamma$ and $f$ which
guarantee the convergence of the integral. Typically we impose
sufficient decay condition on $f$ and some geometric assumption on
$\gamma$. This map $I$ sends a function on $M$ to a function on the
space of geodesics on $M$. It turns out to be useful to instead
consider $I$ as a map from $M$ to the unit sphere bundle $SM$ of $M$
by identifying $\beta\in SM$ with the unique unit speed geodesic whose lift
goes through $\beta$; of course for different $\beta$'s on the same lifted
geodesic $If(\beta)$ is the same.

One of the 
reasons for the importance of this problem is that it is the 
linearization of the boundary rigidity problem, i.e.\ whether
Riemannian metrics can be determined from the (renormalized in many
settings) lengths of their geodesics, see for instance
\cite{uhlmann2016journey}; however, X-ray transforms show up in many
other problems of interest.

The inverse problem concerns the question whether one can determine
$f$ from $If$, i.e., whether $I$ is left invertible, potentially with
additional stability questions (continuity properties of a left
inverse), as well as whether a left inverse can be constructed effectively. The answer
depends on $(M,g)$ and on the function class we choose.
The most famous and `standard' conjecture in this field is Michel's,
namely that boundary rigidity holds on (compact) simple
manifolds. Here a Riemannian manifold with boundary $(M,g)$ is called
simple if for any $p \in M$, the exponential map $\exp_p$ is a
diffeomorphism from a neighborhood of the origin of $T_pM$ and if
$\partial M$ is strictly convex with respect to $g$.

In this paper we consider the geodesic X-ray transform on
asymptotically conic spaces.  Recall that a conic metric, on a manifold
$(0,\infty)_r\times Y$, with $Y$ the cross section or link, which we
always assume is compact and without boundary, is one of the form
\begin{equation}\label{conic_metric_1}
g_\infty=dr^2+r^2 g_0,
\end{equation}
where $g_0$ is a Riemannian metric on $Y$. An asymptotically conic
metric is one on a manifold which outside a compact set is identified
with $(r_0,\infty)_r\times Y$, with a metric that on this conic end
tends to $g_\infty$ as $r\to\infty$ in a specified way. An example is
the Euclidean metric, for which the cross section is the standard
sphere, and indeed metrics asymptotic to the Euclidean one at
infinity, or more generally to the perturbations of the Euclidean
metric by changing the metric on the link, namely the sphere `at
infinity'.

To be
concrete, for our purposes, it is useful to compactify our space,
i.e.\ let $x=r^{-1}$, so $r\to\infty$ corresponds to $x\to 0$, and add a
boundary $\{0\}_x\times Y$ to the manifold, thus compactifying it to
$\overline{M}$. An
asymptotically conic metric then, as introduced by Melrose
\cite{Melrose1994}, is a Riemannian metric on $M$ which is of the form
\begin{equation}\label{conic_metric_2}
g=\frac{dx^2}{x^4}+\frac{\tilde g}{x^2}
\end{equation}
near $\partial\overline{M}$, where $\tilde g$ is a smooth symmetric 2-cotensor on
$\overline{M}$; $g$ is thus asymptotic to $g_\infty$ given by
$\tilde g|_{x=0}$ on the cross section $Y$.

A key difficulty in analyzing the X-ray transform in general is the potentially
complicated geometry, such as the presence of conjugate points, though
these do not exist under Michel's hypotheses. However, on perturbations of
asymptotically Euclidean metrics (for which the link has conjugate
points at distance $\pi$), one typically has
conjugate points; indeed this is necessarily the case if the metric
keeps being asymptotic to Euclidean space but is not flat, as shown
recently by Guillarmou, Mazzucchelli and Tzou
\cite{Guillarmou-Mazzucchelli-Tzou:Conjugate}. Under the assumption of
the absence of conjugate points (as well as other assumptions), Guillarmou, Lassas and Tzou
\cite{Guillarmou-Lassas-Tzou:Conic} have indeed analyzed the the
geodesic X-ray transform on asymptotically conic spaces, but these
results in particular do not apply to (non-trivial) asymptotically Euclidean metrics.

One way of dealing with
the geometric complications, introduced by Uhlmann and Vasy in \cite{uhlmann2016inverse}, is by working locally in smaller regions, namely on super-level
sets of a function whose level sets are strictly concave from the side of the
super-level sets. Such functions always exist locally near the
boundary in Michel's setting due to the strict convexity of the
boundary; indeed the latter guarantees this without further
assumptions. The global existence is called a {\em convex foliation
  condition}, and is satisfied under various conditions discussed in \cite{uhlmann2016inverse,stefanov2018inverting}, see
also a thorough study in
\cite[Section~2]{paternain2019geodesic} and references
therein. The chosen level set acts as an {\em artificial boundary},
analytically pushing to infinity the geometrically finite boundary. We
recall the main result here. For $O$ an open set in $M$, the $O-$local geodesic X-ray transform is the X-ray transform restricted to geodesic segments which are completely in $O$ and in addition with endpoints on $\partial M$.

\begin{thm}[Uhlmann and Vasy, \cite{uhlmann2016inverse}]
Suppose $(M,g)$ is a Riemannian manifold of dimension $\geq 3$ with strictly convex boundary
and $O$ is an open set in it. The $O-$local geodesic X-ray transform
is left invertible on a small collar neighborhood of $\partial M$. It
is globally left invertible under a global convex foliation condition.
\end{thm}

More recently Vasy \cite{vasy2020semiclassical} introduced a
semiclassical approach to this problem, which could also be combined
with the artificial boundary method to take advantage of the best
features of both. Subsequently Vasy and Zachos extended this to
asymptotically conic manifolds in \cite{zachos2022inverting}. In this
extension, it is shown that if we insert a localizer $\tilde{\chi}$ in
the definition of the normal operator (roughly $I^*\tilde\chi I$) and conjugate it by a suitable
exponential weight $e^{\frac{\F\Phi}{h}}$ to define the modified
normal operator, then this operator is an elliptic, and thus for
sufficiently small $h$ invertible, member of a new
pseudodifferential algebra, whose non-semiclassical version, the
1-cusp algebra, already
appeared in Zachos' PhD thesis \cite{Zachos:Thesis}.

\begin{thm}[Vasy and Zachos, \cite{zachos2022inverting}]\label{thm:Zachos}
Suppose that $M$ is a manifold of dimension $n \geq 3$, $g$ is an asymptotically conic metric on $M$ with cross sections without conjugate points within distance $\leq \frac{\pi}{2}$.  
Then on a collar neighborhood of infinity the geodesic X-ray transform is injective on the restriction to the collar neighborhood of functions in $e^{-\frac{C}{x^{2p}}}L_g^2$, where $p>0$ and $C$ depends on $p$.
\end{thm}

Notice that the hypotheses of the theorem in particular allow
asymptotically Euclidean metrics, or more general perturbations at
infinity of asymptotically
Euclidean metrics. However, this comes at the cost of more stringent
decay assumptions on the unknown function $f$ relative to \cite{Guillarmou-Lassas-Tzou:Conic}.

In this paper, we extend this conic result to the one form and 2-tensor
cases. These tensorial problems have an additional complication
relative to the scalar ones above. Potential tensors, i.e.\ tensors of the form $d^s v$
with $v$ a one lower rank tensor vanishing on $\partial M$, or, in our
case, sufficiently fast at infinity, i.e.\ at $\partial\overline{M}$,
where $d^s$ is the symmetrization of the gradient with respect to $g$,
are in the kernel of $I$. In the case of 1-forms, where $d^s$ is the
exterior derivative on functions (so is independent of $g$), this is an immediate
consequence of the
fundamental theorem of calculus. Thus, the natural injectivity we may
expect, which is called {\em s-injectivity}, is that $If=0$ implies
that $f$ is a potential tensor. While this is naturally phrased in
terms of quotient spaces (quotienting by potential tensors),
analytically it is much easier to work in a complementary space to
potential tensors. This space can be obtained by imposing a {\em gauge
  condition}. The `standard' gauge for this problem is the solenoidal
one, namely that $\delta^sf=0$, where $\delta^s$ is the adjoint of the
symmetric gradient with respect to $g$, i.e.\ the (negative) divergence. In the present
situation we need a modification of this gauge condition to one of the
form $\delta^s_{h,\F}f=0$, where $\delta^s_{h,\F}$ is a `Witten type'
divergence, a version of which (adapted to the artificial boundary
there) was introduced in the work of Stefanov, Uhlmann and Vasy
\cite{stefanov2018inverting}, where the scalar local invertibility
result \cite{uhlmann2016inverse} was extended to tensors. Here we use
yet another version to deal with both an artificial boundary
(corresponding to the collar neighborhood) and the asymptotically
conic infinity; see Section~\ref{sec_gauge}.
Our main theorem is:

\begin{thm}
  Suppose $(M,g)$ is as in Theorem~\ref{thm:Zachos}.
Let $\F>0$ for one forms, and $\F$ is sufficiently large for two
tensors, $h>0$ is sufficiently small, and also that $\Omega_{x_0} = \{ x \leq x_0 \}$ with $x_0$ small.
The X-ray transform restricted to geodesics staying in $\Omega_{x_0}$
is injective on restrictions to $\Omega_{x_0}$ of tensors decaying
with rate $e^{-\frac{C}{x^{2}}}$ with $C=\F/(2h)$ (in the sense of membership
in $e^{-\frac{C}{x^2}}L^2$) and satisfying the gauge condition
 \begin{align}
\delta_{h,\F}^s (e^{-\frac{\F\Phi}{h}}f) = 0,
\label{gauge_eq}
 \end{align}
 where $\delta^s_{h,\F}$ is given in \eqref{d_delta_defn} and $\Phi$
 in \eqref{Phi_defn}.
 \label{thm_conic}
\end{thm}

\begin{rmk}
  Much as for Theorem~\ref{thm:Zachos}, the result continues to hold when the decay requirement is weakened to
be the rate $e^{-\frac{C}{x^{2p}}}$ for any $p>0$. We only give a detailed proof for the case $p=1$ as stated above, and indicate minor changes needed for the general $p$ case in
Remark \ref{rmk_general_p}.
\label{rmk_general_p_0}
\end{rmk}

Next, let us recall the origin of $\pi/2$ in the last two stated
theorems. A computation of Melrose and Zworski \cite{RBMZw} shows that for an actual
conic metric, such as $g_\infty$, unit speed (prior to reparameterization) geodesics can be
explicitly described explicitly as follows. Here the description will
be the Hamiltonian one, i.e.\ using the cosphere bundle rather than
the sphere bundle. Writing covectors as
$$
-\tau\,dr+\mu\cdot (r\,dy)=\tau\frac{dx}{x^2}+\mu\cdot\frac{dy}{x},
$$
$y$ local coordinates on the link $Y$,
thus in a way adapted to the asymptotically conic structure (these
are coordinates on the scattering cotangent bundle in Melrose's terminology),
bicharacteristics $\gamma$ can be written
as follows using the notation $\mu=|\mu|_{g_0}\hat\mu$:
\begin{equation}\begin{aligned}\label{eq:conic-bichar}
    &    x=\frac{x_0}{\sin r_0}\sin (r+r_0),\ \tau=\cos (r+r_0),\ |\mu|=\sin (r+r_0),\\
    & (y,\hat\mu)=\exp(rH_{\frac{1}{2} h})(y_0,\hat\mu_0),\ r\in(-r_0,-r_0+\pi),
  \end{aligned}\end{equation}
with $(y,\hat\mu)$ thus following a unit speed lifted geodesic of
length $\pi$ in
$Y$. Note that the maximum of $x\circ \gamma$, which is the point of tangency to
level sets of the function $x$, occurs halfway in the domain of
$\gamma$, at (parameter) distance $\pi/2$ from either endpoint, at $r+r_0=\pi/2$, and
thus in terms of the boundary geodesic distance $\pi/2$ from either
endpoint. Due to our exponential weights and cutoffs, the key point is
to have no points conjugate to the point of tangency to this level set
along these geodesics, which is guaranteed if the link has no
conjugate points within distance $\pi/2$. While this is for actual
conic metrics, for asymptotically conic metrics the analogous
condition automatically holds in a sufficiently small collar
neighborhood of the boundary.

The condition \eqref{gauge_eq}
is not restrictive in the sense that
we can arrange it by adding
a potential tensor, which does not affect the result of the X-ray
transform (i.e.\ is in its kernel); in this sense (for suitably fast
decaying tensors) Theorem~\ref{thm_conic} is optimal.

\begin{thm}\label{thm:gauging}
Suppose that $\F>0$ is sufficiently large and $h>0$ is sufficiently small,
and let $C=\F/(2h)$.
For each one form and 2-tensor $f$ decaying with rate
$e^{-\frac{C}{x^2}}$ (as $x \rightarrow 0$, in the sense of membership
in $e^{-\frac{C}{x^2}}L^2$), there exists a tensor $v$ of one lower rank and the same
exponential decay rate such that 
\begin{align*}
\delta^s_{h,\F}(e^{-\frac{\F\Phi}{h}} (f-d^sv)) = 0.
\end{align*}
\label{thm_gauge}
\end{thm}

The structure of this paper is as follows. In Section~\ref{sec:psdos}
we recall the 1-cusp and scattering pseudodifferential algebras from
the scalar setting of \cite{zachos2022inverting}; here going to
tensors does not cause any complications. Note that the algebras used
do {\em not} match the geometry in the sense that for instance an
asymptotically conic metric naturally corresponds to Melrose's
scattering algebra \cite{Melrose1994} (e.g.\ its Laplacian is in this
class), but we instead use a different, 1-cusp, algebra
there. Ultimately the reason we {\em can} do this is that the X-ray
transform problem is overdetermined if $n\geq 3$: {\em the way we choose
what information to keep via the cutoff $\tilde\chi$} determines the
operator algebra structure in tandem with the geometry.
In
Section~\ref{sec:geometry} we thus analyze the geometric operators such as
the symmetric gradient and divergence as elements of the new
algebras. In Section~\ref{sec:normal} we analyze the modified, gauge fixed, normal
operator, showing that it is an elliptic element of the algebra, and
thus for sufficiently small semiclassical parameter $h$ it is
invertible, thus proving the main theorem,
Theorem~\ref{thm_conic}. Finally in
Section~\ref{sec_gauge} we show that the gauge condition can be
arranged, thus proving Theorem~\ref{thm:gauging}. Note that arranging the gauge
condition in the present setting is {\em much} easier than in
\cite{stefanov2018inverting}, since in the latter paper the {\em actual}
boundary of the manifold (as opposed to the artificial boundary)
caused significant complications.

\section{The 1-cusp and scattering pseudodifferential algebras}\label{sec:psdos}
In this section we recall the analytic ingredients, namely the relevant
pseudodifferential algebras, from \cite{zachos2022inverting}.

\subsection{The semiclassical foliation 1-cusp algebra}
We briefly describe the 1-cusp pseudodifferential algebra here. For
detailed construction and explanation, see Section 2 of
\cite{zachos2022inverting}. First we define the cusp vector fields and
the 1-cusp vector fields. Let $x$ be a boundary defining function of
$M$. Then $\mathcal{V}_{\mathrm{cu}}$ consists of smooth vector fields
$V$ tangent to $\partial M$ such that $Vx = O(x^2)$. Note that the
this is a Lie algebra of vector fields (under commutators) and it
depends on the choice of $x$ modulo $O(x^2)$. In local coordinates $x,y_1,...,y_{n-1}$, they have the form
\begin{align*}
a_0(x,y)x^2D_x + \sum_{j=1}^{n-1} a_j(x,y)D_{y_j},
\end{align*}
where $a_j$ are smooth functions of their variables. Then we define the 1-cusp vector fields as cusp vector fields with one extra vanishing order near the boundary:
\begin{align*}
\mathcal{V}_{\mathrm{1c}}(M): = x\mathcal{V}_{\mathrm{sc}}(M).
\end{align*}
In local coordinates, they have the form
\begin{align*}
a_0(x,y)x^3D_x + \sum_{j=1}^{n-1} a_j(x,y)xD_{y_j}.
\end{align*}
Over $C^\infty(M)$, this generates the algebra of 1-cusp differential
operators as (locally) finite sums of finite products of these.

Next we introduce the semiclassical foliation algebra associated to
$\mathcal{F}$, the foliation given by the level sets of $x$; this
depends on $x$ even more strongly since it depends on all of the level
sets of $x$. The usual the semiclassical version of
$\mathcal{V}_{\mathrm{1c}}(M)$ is defined to be
$\mathcal{V}_{\mathrm{1c},h}(M)=h\mathcal{V}_{\mathrm{1c}}(M)$. Its
variant, the Lie algebra of semiclassical foliation vector fields is
\begin{align*}
\mathcal{V}_{\mathrm{1c},h,\mathcal{F}}(M) = h\mathcal{V}_{\mathrm{1c}}(M) + h^{1/2}\mathcal{V}_{\mathrm{1c}}(M;\mathcal{F}),
\end{align*}
where $\mathcal{V}_{\mathrm{1c}}(M;\mathcal{F})$ is the collection of
1-cusp vector fields tangent to the foliation, i.e.\ locally of the
form
\begin{align*}
\sum_{j=1}^{n-1} a_j(x,y)xD_{y_j}.
\end{align*}
This again generates a differential operator algebra;
a typical semiclassical foliation 1-cusp differential operator is thus of the form
\begin{align*}
\sum_{\alpha+|\beta| \leq m} a_{\alpha\beta}(x,y,h)(hx^3D_x)^\alpha(h^{1/2}xD_y)^\beta.
\end{align*}
Correspondingly, the frame of the semiclassical foliation 1-cusp cotangent bundle, which is denoted by 
$_{h,\mathcal{F}}^{\;\;\mathrm{1c}}T^*X$ (and its scattering counterpart is denoted by $_{h,\mathcal{F}}^{\;\;\mathrm{sc}}T^*X$), is given by
\begin{align*}
\frac{dx}{hx^3},\frac{dy_j}{h^{1/2}x}.
\end{align*}
Suppose coordinates of this bundle in this frame are written as $\xi_{\mathrm{1c}},\eta_{\mathrm{1c}}$, then our symbol class $S^{m,l}_{\mathrm{1c},h,\mathcal{F}}(M)$ consists of standard semiclassical symbols in this coordinate system:
\begin{align*}
|(xD_x)^\alpha D_y^\beta D_{\xi_{\mathrm{1c}}}^\gamma D_{\eta_{\mathrm{1c}}}^\delta a(x,y,\xi_{\mathrm{1c}},\eta_{\mathrm{1c}},h)| \leq C_{\alpha\beta\gamma\delta} \la (\xi_{\mathrm{1c}},\eta_{\mathrm{1c}}) \ra^{m-\gamma-\delta}x^{-l}.
\end{align*}
Although in coordinates this is the same definition as the cusp symbol class, in fact they are symbols on a different cotangent bundle, and correspondingly they are quantized in a different manner. The quantization map is defined by
\begin{align}
A_hu(x,y) = & (2\pi)^{-n}h^{-n/2-1/2}
\\& \int e^{i (\frac{x-x'}{x^3} \frac{\xi_{\mathrm{1c}}}{h}+ \frac{y-y'}{x} \frac{\eta_{\mathrm{1c}}}{h^{1/2}})} a(x,y,\xi_{\mathrm{1c}},\eta_{\mathrm{1c}},h) u(x',y') \frac{dx'dy'}{(x')^{n+2}}d\xi_{\mathrm{1c}}d\eta_{\mathrm{1c}};
\label{1c_quant_defn}
\end{align}
this is understood to be valid
away from $\{x=0\}$, this is identical to the `standard' semiclassical
foliation operators intorduced in \cite{vasy2020semiclassical}, while
for $h>0$, this gives the 1-cusp pseudodifferential operators
introduced by Zachos \cite{Zachos:Thesis}. The behavior near $x=h=0$
is the important point for us.

Next we define the ellipticity of symbols and operators.
\begin{defn}
A symbol $a \in S^{m,l}_{\mathrm{1c},h,\mathcal{F}}(M)$ is called elliptic if
\begin{align*}
|a(x,y,\xi_{\mathrm{1c}},\eta_{\mathrm{1c}})| \geq cx^{-l} \la (\xi_{\mathrm{1c}},\eta_{\mathrm{1c}}) \ra^m, \quad c>0;
\end{align*}
its quantization $A$ is also called elliptic in this case.
\end{defn}
Under this condition, see \cite[Section~2.5]{zachos2022inverting}, its quantization $A$ has a parametrix $B \in \Psi_{\mathrm{1c},h,\mathcal{F}}^{-m,-l}(M)$ such that
\begin{align*}
AB-\Id, BA-\Id \in h^{\infty}\Psi_{\mathrm{1c},h,\mathcal{F}}^{-\infty,-\infty}(M).
\end{align*}

One can now define the semiclassical foliation 1-cusp Sobolev spaces
$H^{s,r}_{\mathrm{1c},h,\mathcal{F}}(M)$, see
\cite[Section~2.5]{zachos2022inverting}, for instance for $s\geq 0$ by choosing
$A\in \Psi_{\mathrm{1c},h,\mathcal{F}}^{s,0}(M)$ elliptic, and demanding
$$
u\in H^{s,r}_{\mathrm{1c},h,\mathcal{F}}(M) \Leftrightarrow u\in x^r L^2(M)\ \text{and}\ Au\in x^r L^2(M);
$$
here $L^2(M)$ is the $L^2$ space relative to a fixed polynomially
weighted density, which in the geometric context of asymptotically
conic spaces is natural to take to be the metric density, which is
equivalent to $\frac{dx\,dy}{x^{n+1}}$. Equivalently, for $s\geq 0$
integer,
$$
\|u\|^2_{H^{s,r}_{\mathrm{1c},h}}=\|x^{-r}u\|^2+\sum_{j+|\alpha|\leq s}\|(hx^3D_x)^j(h^{1/2}xD_y)^\alpha\|^2
$$
with the spaces for other $s$ defined via interpolation and duality.

Pseudodifferential operators are bounded on
these Sobolev spaces, namely for all $s,r$,
$$
A \in \Psi_{\mathrm{1c},h,\mathcal{F}}^{m,l}(M)\Rightarrow A\in
\mathcal{L}(H^{s,r}_{\mathrm{1c},h,\mathcal{F}}, H^{s-m,r-l}_{\mathrm{1c},h,\mathcal{F}}).
$$
The existence of parametrices for elliptic operators implies that for elliptic operarors
$A$, there exists $h_0>0$ such that for $h \in [0,h_0)$, $A\in
\mathcal{L}(H^{s,r}_{\mathrm{1c},h,\mathcal{F}}, H^{s-m,r-l}_{\mathrm{1c},h,\mathcal{F}})$ is
invertible with uniform bounds. This is a key reason to use
the semiclassical algebra: the errors are not only compact, but can be
indeed eliminated via a convergent Neumann series.

\subsection{The semiclassical foliation scattering algebra}
Now we recall basic facts about the semiclassical foliation scattering
algebra from \cite{vasy2020semiclassical}. This pseudodifferential
operator algebra is similar to its 1-cusp analogue. Let
$\mathcal{V}_{\mathrm{sc}}(M)$ be scattering vector fields, i.e.\ $x$
times vector fields tangent to the boundary, so in local coordinates
elements are of the form
\begin{align*}
a_0(x,y)x^2D_x + \sum_{j=1}^{n-1} a_j(x,y)xD_{y_j},
\end{align*}
where $a_j$ are smooth functions of their variables.
Also let
$\mathcal{V}_{\mathrm{sc}}(M;\mathcal{F})$ be those ones tangent to
level sets of the foliation, i.e. of the form
\begin{align*}
\sum_{j=1}^{n-1} a_j(x,y)xD_{y_j},
\end{align*}
The vector fields of relevance are combinations of $h-$semiclassical or $h^{1/2}-$semiclassical and tangent to the foliation:
\begin{align*}
\mathcal{V}_{\mathrm{sc},h,\mathcal{F}}(M) = h \mathcal{V}_{\mathrm{sc}}(M)+h^{1/2}\mathcal{V}_{\mathrm{sc}}(M;\mathcal{F}).
\end{align*}
Denote the coordinate of fiber part in the frame $\frac{dx}{hx^2},\frac{dy_j}{h^{1/2}x}$ by $\xi_{\mathrm{sc}},\eta_{\mathrm{sc}}$. The symbol class $S^{m,l}_{\mathrm{sc},h,\mathcal{F}}(M)$ consists of smooth functions such that
\begin{align*}
|(xD_x)^\alpha D_y^\beta D_{\xi_\mathrm{sc}}^\gamma D_{\eta_\mathrm{sc}}^\delta a(x,y,\xi_\mathrm{sc},\eta_\mathrm{sc},h)| \leq C_{\alpha \beta \gamma \delta} \la (\xi_\mathrm{sc},\eta_\mathrm{sc}) \ra^m x^{-l}.
\end{align*}
For such $a$, we quantize it to be $A \in \Psi^{m,l}_{\mathrm{sc},h,\mathcal{F}}(M)$ by
\begin{align}
A_h u(x,y) = & (2\pi)^{-n}h^{-n/2-1/2} 
\\ & \int e^{i ( \frac{x-x'}{x^2} \frac{\xi_{\mathrm{sc}}}{h} + \frac{y-y'}{x} \frac{\eta_{\mathrm{sc}}}{h^{1/2}} ) } a(x,y,\xi_{\mathrm{sc}},\eta_{\mathrm{sc}},h)u(x',y') \frac{dx'dy'}{(x')^{n+1}}d\xi_{\mathrm{sc}} d\eta_{\mathrm{sc}}.
\label{sc_quant_defn}
\end{align}
Away from $\{x=0\}$, these are just the standard semiclassical foliation operators defined in \cite{vasy2020semiclassical}, and in $h>0$ are the standard scattering pseudodifferential operators, with combined behavior near $x=h=0$. The ellipticity condition changes correspondingly.
\begin{defn}
A symbol $a \in S^{m,l}_{\mathrm{sc},h,\mathcal{F}}(M)$ is said to be elliptic if
\begin{align*}
|a(x,y,\xi_{\mathrm{sc}},\eta_{\mathrm{sc}},h)| \geq cx^{-l} \la (\xi_{\mathrm{sc}},\eta_{\mathrm{sc}})\ra^{m};
\end{align*}
its quantization $A$ is also called elliptic in this case.
\end{defn}
Similar to the 1-cusp case, its quantization $A$ has a parametrix $B \in \Psi_{\mathrm{1c},h,\mathcal{F}}^{-m,-l}(M)$ such that

\begin{align*}
AB-\Id, BA-\Id \in h^{\infty}\Psi_{\mathrm{sc},h,\mathcal{F}}^{-\infty,-\infty}(M),
\end{align*}
and there exists $h_0>0$ such that for $h \in [0,h_0)$, $A \in \mathcal{L}(H^{s,r}_{\mathrm{sc},h,\mathcal{F}}, H^{s-m,r-l}_{\mathrm{sc},h,\mathcal{F}})$ is invertible with uniform bounds. The same as the 1-cusp case, the errors are not only compact, but can be indeed eliminated.

\subsection{The combined class}
\label{the_combined_class}
Our analysis in the asymptotically conic setting is a combination of the 1-cusp (due to the conic end) and the scattering (due to the artificial boundary) structures. Let $X$ be a manifold with boundary equiped with a function $x$ such that $dx$ is never degenerate and level sets $\Sigma_c:=x^{-1}(c)$ are smooth hypersurfaces, among which $\Sigma_0$ and $\Sigma_{x_0}$ are two boundary surfaces. In our application, $\Sigma_0$ will be the infinity of our asymptotic conic manifold, $\Sigma_{x_0}$ will be taken to be the artificial boundary used to localize our argument and $X$ is a domain in $M$.

We define the operator class $\Psi_{\mathrm{sc},\mathrm{1c},h,\mathcal{F}}^{m,l_1,l_2}(X)$ first, which consists of pseudodifferential operators of $m-$th diffrential order that are scattering of order $l_1$ in $\Omega_{x_0}$, which is a neighborhood of $\Sigma_{x_0}$, and are 1-cusp of order $l_2$ in $\Omega_0$, which is a neighborhood of $\Sigma_0$.     

\begin{defn}
The operator class $\Psi_{\mathrm{sc},\mathrm{1c},h,\mathcal{F}}^{m,l_1,l_2}(X)$ consists of operators $A$ such that 
\begin{itemize}
	\item For $\phi,\psi \in  C^{\infty}(X)$ with support disjoint from $\Sigma_0$, $\phi A \psi \in \Psi_{\mathrm{sc},h,\mathcal{F}}^{m,l_1}(\Omega_{x_0})$, and the semiclassical foliation scattering algebra is constructed using $\Sigma_{x_0}$ as the boundary surface.
	\item For $\phi,\psi \in  C^{\infty}(X)$ with support disjoint from $\Sigma_{x_0}$, $\phi A \psi \in \Psi_{\mathrm{1c},h,\mathcal{F}}^{m,l_2}(\Omega_{0})$, and the semiclassical foliation 1-cusp algebra is constructed using $\Sigma_{0}$ as the boundary surface.
	\item For $\phi,\psi \in  C^{\infty}(X)$ with
          disjoint support, $\phi A \psi$ has Schwartz kernel which is
          $C^\infty$ and rapidly vanishing in the semiclassical
          parameter $h$ as well as at all boundary hypersurfaces of $\bar{\Omega}\times\bar{\Omega}$.
\end{itemize}
\label{defn_operator}
\end{defn}
One can easily check that
$\Psi_{\mathrm{sc},\mathrm{1c},h,\mathcal{F}}^{\infty,\infty,\infty}(X)$ is a
tri-filtered $*-$algebra. We also need Sobolev spaces; as usual these
are defined by localization:
\begin{defn}
The function class $H_{\mathrm{sc},\mathrm{1c},h,\mathcal{F}}^{m,l_1,l_2}(X)$ consists of functions $f$ such that 
\begin{itemize}
	\item For $\phi \in  C^{\infty}(X)$ with support
          disjoint from $\Sigma_0$, $\phi f  \in
          H_{\mathrm{sc},h,\mathcal{F}}^{m,l_1}(\Omega_{x_0})$, and the semiclassical foliation 
          scattering Sobolev space is constructed using $\Sigma_{x_0}$
          as the boundary surface.
	\item For $\phi \in  C^{\infty}(X)$ with support
          disjoint from $\Sigma_{x_0}$, $\phi f  \in
          H_{\mathrm{1c},h,\mathcal{F}}^{m,l_2}(\Omega_{0})$, and the semiclassical foliation 
          1-cusp Sobolev space is constructed using $\Sigma_{x_0}$ as
          the boundary surface.
\end{itemize}
\label{defn_function}
\end{defn}

We also define a new cotengent bundle that has corresponding boundary behavior near each boundary surface.
\begin{defn}
The semiclassical sc-1c foliation cotangent bundle $_{h,\mathcal{F}}^{\mathrm{sc},\mathrm{1c}}T^*X$ is the vector bundle such that
\begin{itemize}
\item In a neighborhood of $\Sigma_{x_0}$, its local frame is given by the frame of the semiclassical foliation scattering cotangent bundle: $\frac{dx}{h(x_0-x)^2}, \frac{dy_j}{h^{1/2}(x_0-x)}$.
\item In a neighborhood of $\Sigma_0$, its local frame is given by the frame of the semiclassical foliation 1-cusp cotangent bundle:  $\frac{dx}{hx^3}, \frac{dy_j}{h^{1/2}x}$.
\item Away from both $\Sigma_0,\Sigma_{x_0}$, its local frame is the
  same as the foliation cotangent bundle: $\frac{dx}{h},\frac{dy_j}{h^{1/2}}$.
\end{itemize}
\label{defn_bundle}
\end{defn}

Bundle valued Sobolev spaces
$H_{\mathrm{sc},\mathrm{1c},h,\mathcal{F}}^{s,l_1,l_2}(X,_{h,\mathcal{F}}^{\mathrm{sc},\mathrm{1c}}T^*X)$
are defined to be spaces of sections of $_{h,\mathcal{F}}^{\mathrm{sc},1c}T^*X$ with
coefficients in $H_{\mathrm{sc},\mathrm{1c},h,\mathcal{F}}^{s,l_1,l_2}(X)$;
similarly for $\mathrm{Sym}_{h,\mathcal{F}}^{2,\mathrm{sc},\mathrm{1c}}T^*X$, which is
the symmetric part of
$_{h,\mathcal{F}}^{\mathrm{sc},1c}T^*X\otimes\leftidx{_{h,\mathcal{F}}^{\mathrm{sc},\mathrm{1c}}}{T^*X}$. Moreover,
the operator classes mapping between those bundles, acting as a linear
map with components in
$\Psi_{\mathrm{sc},\mathrm{1c},h,\mathcal{F}}^{m,l_1,l_2}(X)$, are denoted by
\begin{align*}
\Psi_{\mathrm{sc},\mathrm{1c},h,\mathcal{F}}^{m,l_1,l_2}(X;_{h,\mathcal{F}}^{\mathrm{sc},\mathrm{1c}}T^*X,_{h,\mathcal{F}}^{\mathrm{sc},\mathrm{1c}}T^*X)
\end{align*}
and
\begin{align*}
\Psi_{\mathrm{sc},\mathrm{1c},h,\mathcal{F}}^{m,l_1,l_2}(X;\mathrm{Sym}_{h,\mathcal{F}}^{2,\mathrm{sc},\mathrm{1c}}T^*X,\mathrm{Sym}_{h,\mathcal{F}}^{2,\mathrm{sc},\mathrm{1c}}T^*X)
\end{align*}
respectively.

\section{The geometric operators as elements of the pseudodifferential algebras}\label{sec:geometry}
In this section we analyze the operators in the gauge condition. Let
$\nabla$ be the covariant derivative with respect to $g$, $d^s$ be the
symmetrization of $\nabla$; these are natural geometric objects that
are important to keep unchanged since $d^s$ gives rise to the kernel
of the X-ray transform. On the other hand, the operators involved in
the gauge condition are artificial, and we need to (and can) choose
them to match the analytic framework. Keeping in mind that we employ
an analytic framework which is 1-cusp at the conic infinity and
scattering at the artificial boundary, we choose
$g_{\mathrm{sc},\mathrm{1c},h}$ to be a combination 
of the semiclassical scattering metric and the semiclassical 1-cusp metric,
and let $\delta^s$ be the adjoint of
$d^s$ with respect to $g_{\mathrm{sc},\mathrm{1c},h}$, so it is the (negative) divergence operator.
Concretely in a neighborhood of $\Sigma_0$, we have
\begin{align}
g_{\mathrm{sc},\mathrm{1c},h} = h^{-2}x^{-6}dx^2 + h^{-1}x^{-2}g_1, \label{metric_sc_1c_1}
\end{align}
while  in a neighborhood of $\Sigma_{x_0}$, we have
\begin{align}
g_{\mathrm{sc},\mathrm{1c},h} = h^{-2}(x_0-x)^{-4}dx^2+h^{-1}(x_0-x)^{-2}g_2. \label{metric_sc_1c_2}
\end{align}
$g_1,g_2$ are smooth families of Riemannian metric on level sets of $x$ in neighborhoods mentioned above.
The specific choice among smooth transitions between those two
boundary faces does not affect our analysis and we choose one of them
and fix it. We reiterate that the metric $g_{\mathrm{sc},\mathrm{1c},h}$ is introduced as an analytic tool to define adjoint operators,
 convert tensors and combine the analysis near $\Sigma_0$ and $\Sigma_{x_0}$.
In particular, near $\Sigma_0$, $\delta^s$ is the adjoint of $d^s$ with respect to the semiclassical foliation 1-cusp metric $h^{-2}x^{-6}dx^2 + h^{-1}x^{-2}g_1$, and near $\Sigma_{x_0}$  $\delta^s$ is the adjoint of $d^s$ with respect to the semiclasscial foliation scattering metric $h^{-2}(x_0-x)^{-4}dx^2+h^{-1}(x_0-x)^{-2}g_2$.

Next we conjugate them by an exponential weight so that they have desired analytic properties.
First we define
\begin{align}
\Phi = F_0 \circ x,  \label{Phi_defn}
\end{align}
where $F_0$ is a smooth, increasing function, in the strong sense that $F'_0>0$, such that $F_0(x) = -\frac{1}{2x^2}$ near $\Sigma_0$ and $F_0 = \frac{1}{x_0-x}$ near $\Sigma_{x_0}$. Then our conjugated symmetric differential and its adjoint are
\begin{align}
d^s_{h,\F} = e^{-\frac{\F\Phi}{h}} d^s e^{\frac{\F\Phi}{h}}, \quad  \delta^s_{h,\F} = e^{\frac{\F\Phi}{h}} \delta^s e^{-\frac{\F\Phi}{h}}.
\label{d_delta_defn}
\end{align}
We then consider the effect of conjugation when we compute symbols in the 1-cusp algebra. Let $\Phi(x)=-\frac{1}{2x^2}$, we have                    
\begin{align}
\begin{split}
e^{-\frac{F\Phi(x)}{h}}hx^3D_xe^{ \frac{F\Phi(x)}{h}} &= e^{\frac{F}{2hx^2}}hD_xe^{-\frac{F}{2hx^2}}\\
 & = hx^3D_x-iF.
 \end{split} \label{conjugation_conic}
\end{align} 
Thus the effect of exponential conjugation is replacing $\xi$ by $\xi-i\F$. We then write the exterior derivative $d_0$ in terms of semiclassical foliation 1-cusp covectors:
\begin{align*}
d_0f = (\partial_xf) dx +  \sum_j (\partial_{y_j}f)dy_j = (hx^3\partial_xf) \frac{dx}{hx^3} + \sum_j (h^{1/2}x\partial_{y_j})\frac{dy_j}{h^{1/2}x}.
\end{align*}
This already shows that $d_0$, which coincides with $\nabla$ when acting on functions, has principal symbol 
$\xi_{\mathrm{1c}}\frac{dx}{hx^3}\otimes+\eta\cdot\frac{dy}{h^{1/2}x}\otimes$ when considered as a first order semiclassical foliation 1-cusp differential operator. After conjugation, as mentioned after (\ref{conjugation_conic}), $\xi_{\mathrm{1c}}$ is replaced by $\xi_{\mathrm{1c}}-i\F$. A similar computation shows that its principal symbol when viewed as operator on tensors with higher rank is also tensoring with the covector at which the symbol is evaluated.

\subsection{$d^s$ as a scattering differential operator}
\label{section_Chris} 
Next we consider the action of $d^s$ on one forms, which originally is a first order differential operator sending sections of $\leftidx{^{\mathrm{sc}}}{T^*X}$ to sections of $\mathrm{Sym}^{2,\mathrm{sc}}T^*X$. We compute its principal symbol in $\Psi_{\mathrm{sc}}^{1,0}$. We consider $\nabla$ first, whose action on a  scattering one form $T=\hat{T}_0\frac{dx}{x^2}+\sum_{j=1}^{n-1} \hat{T}_j \frac{dx^j}{x}$ is given by
\begin{align}
\nabla(T) = \partial_cT_b dx^c \otimes dx^b - \Gamma^d_{bc}T_d dx^c \otimes dx^b,  \label{nabla}
\end{align}
where we used the Einstein's convention of summation, $T_0=x^{-2}\hat{T}_0, T_j = x^{-1}\hat{T}_j$ for $1 \leq j \leq n-1$, and $\Gamma^d_{bc}$ is the Christoffel symbol defined by
\begin{align}
\Gamma^l_{jk} = \frac{1}{2} g^{lr}( \partial_kg_{rj}+\partial_jg_{rk}-\partial_rg_{jk} ), \label{christoffel_symbol}
\end{align}
where $g^{lr}$ are components of the dual metric. In order to compute
the form of $\Gamma^l_{jk}$, we consider the dual metric first. Recall (\ref{conic_metric_2}), written as a block matrix, with respect to the scattering basis, $g_{ij}$ has the form
\begin{align}
\begin{pmatrix}
1   && O(x) \\
O(x) && \tilde{g}
\end{pmatrix},
\label{metric_matrix}
\end{align}
where two $O(x)$ blocks of the shape $1 \times n$ and $n \times 1$
respectively are because the $dx\otimes dy$ terms encoded in
$x^{-2}\tilde{g}$ in (\ref{conic_metric_2}) are of order $O(x)$ when
written in terms of $\frac{dx}{x^2}\otimes \frac{dy}{x}$, and $O(x)$
stands for $x$ times a matrix of the appropriate type with smooth entries.
Inverting this matrix, the dual metric tensor in terms of $x^2\partial_x \otimes x^2\partial_x, x^2\partial_x\otimes \partial_y, x\partial_y \otimes x\partial_y$ has the form
\begin{align}
\begin{pmatrix}
1   && O(x) \\
O(x) && \tilde{g}^{-1}
\end{pmatrix}.
\label{dualmetric_matrix}
\end{align}
Thus, in terms of the basis given by tensor products of $\partial_x,\partial_{y_i}$, we know
\begin{align}
\begin{split}
& g^{00}=x^4, \quad g^{0i}=x^4\tilde{g}^{0i}, \quad g^{i0}=x^4\tilde{g}^{i0},  1 \leq i \leq n-1,\\
& g^{ij}=x^2\tilde{g}^{ij} , \quad 1 \leq i,j \leq n-1,     
\end{split}                \label{dual_metric}
\end{align} 
where $\tilde{g}^{0i},\tilde{g}^{i0},\tilde{g}^{ij}$ are functions
smooth down to $x=0$. We consider (\ref{nabla}), broken up into
several cases.

Case 1: $c=0,b=0$ \\
In this case, using (\ref{christoffel_symbol}),  we have (when $b,c$
have fixed value, the summation convention does {\em not} apply to them)
\begin{align*}
& \partial_cT_bdx^c \otimes dx^b 
\\ = & \partial_0T_0 dx^0\otimes dx^0 
\\ = & \partial_x (x^{-2}\hat{T}_0 )dx^0\otimes dx^0 
\\ = & (x^{-2}\partial_x\hat{T}_0-2x^{-3}\hat{T}_0)dx^0\otimes dx^0 
\\ = & x^2\partial_x\hat{T}_0 \frac{dx^0}{x^2} \otimes \frac{dx^0}{x^2}-x\hat{T}_0\frac{dx^0}{x^2} \otimes \frac{dx^0}{x^2},
\\& \Gamma^d_{bc}T_d = \Gamma^d_{00}T_d,
\\& \Gamma^d_{00} = \frac{1}{2} g^{dr}( \partial_xg_{r0}+\partial_xg_{r0}-\partial_rg_{00} ).
\end{align*}
We compute the power of $h$ and $x$ of $\Gamma^d_{00}T_ddx^0 \otimes dx^0 $, when written as an semiclassical foliation 1-cusp tensor:\\
The term $d=0$ is
 $\Gamma^0_{00}h^{-1}x^{-2}\hat{T}_0dx^0 \otimes dx^0=x^2\Gamma^0_{00}\hat{T}_0\frac{dx}{x^2} \otimes \frac{dx}{x^2}$. In the expression of $\Gamma^0_{00}$, the term with $r=0$ contributes $x^4\tilde{g}^{d0}\partial_xx^{-4}=O(x^{-1})$, while terms with $r \neq 0$ contribute $x^4\tilde{g}^{0r}\partial_x(x^{-2}\tilde{g}_{r0})=O(x)$. Thus the $r=0$ case gives the main contribution and combining the $x^2$ factor in the front, $d=0$ term is of order $O(x)$.

When $d \neq 0$, we have $\Gamma^d_{00}T_ddx^0\otimes dx^0 = h^{3/2}x^3\Gamma^d_{00}\hat{T}_d\frac{dx}{x^2}\otimes \frac{dx}{x^2}$. For $\Gamma^d_{00}$, the term $r=0$ is $\frac{1}{2}x^4\tilde{g}^{d0}\partial_xx^{-4}=O(x^{-1})$, while terms with $r \neq 0$ give $x^2\tilde{g}^{dr}\partial_x(x^{-2}\tilde{g}_{r0}) = O(x^{-1})$. Combining with $h^{3/2}x^3$ in the front, it is $O(h^{3/2}x^2)$.

Combining $d \neq 0$ and $d=0$ cases, the coefficient of $\frac{dx}{x^2}\otimes \frac{dx}{x^2}$, together with the term $x\hat{T}_0\frac{dx^0}{x^2} \otimes \frac{dx^0}{x^2}$ created by commuting the $x$ factor and differentiation, is $O(x)$.

Case 2: $c=0,b\neq 0$.
\begin{align*}
& \partial_0T_bdx^0 \otimes dx^b 
\\ = & (\partial_x x^{-1}\hat{T}_b) dx^0 \otimes dx^b
\\ = & x^{-1}\partial_x\hat{T}_b dx^0 \otimes dx^b -x^{-2}\hat{T}_b dx^0\otimes dx^b
\\ = & (x^2\partial_x \hat{T}_b) \frac{dx^0}{x^2}\otimes \frac{dx^b}{x}
-x\hat{T}_b  \frac{dx^0}{x^2}\otimes \frac{dx^b}{x},\\
& \Gamma^d_{bc}T_d = \Gamma^d_{b0}T_d,\\
& \Gamma^d_{b0} = \frac{1}{2} g^{dr}( \partial_xg_{rb}+\partial_bg_{r0}-\partial_rg_{b0} ).
\end{align*}
We compute the power of $h$ and $x$ of $\Gamma^d_{b0}T_ddx^0 \otimes dx^b $, when written as a scattering semiclassical foliation tensor:\\

Consider $d=0$ and $d \neq 0$ respectively. The term with $d=0$ is
\begin{align*}
\Gamma^0_{b0}T_0dx^0 \otimes dx^b = x\Gamma^0_{b0} \hat{T}_0 \frac{dx}{x^2} \otimes \frac{dx^b}{x},
\end{align*}
For terms in the bracket in the expression of $\Gamma^0_{b0}$, the only possible term of $O(x^{-4})$ is $\partial_bg_{r0}$ with $r=0$, which however vanishes since $g_{00}=x^{-4},b \neq 0$. So the term in the bracket is at most of order $O(x^{-3})$. Combining with $g^{0r}=O(x^4)$, we know $\Gamma^0_{b0}=O(x)$, thus $x\Gamma^0_{b0} \hat{T}_0 \frac{dx}{x^2} \otimes \frac{dx^b}{x}=O(x^2)\frac{dx}{x^2} \otimes \frac{dx^b}{x}$.

For $d \neq 0$, we have
\begin{align*}
\Gamma^d_{b0}T_ddx^0 \otimes dx^b = x^2\Gamma^d_{b0} \hat{T}_d \frac{dx}{x^2} \otimes \frac{dx^b}{x}.
\end{align*}
For terms in the expression of $\Gamma^d_{b0}$, as explained above, terms in the bracket are at most $O(x^{-3})$, whereas the dual metric factor has at least $O(x^2)$ vanishing, thus in total we have $\Gamma^d_{b0}=O(x^{-1})$. With the $x^2$ factor in the front, this part has $O(x)$ contribution.
Combining with the case $d=0$, we know the zeroth differential order terms relevant to $\frac{dx}{x^2} \otimes \frac{dx^b}{x}$ component has $O(x)$ scale.

Case 3: $c \neq 0, b=0$.
\begin{align*}
& \partial_cT_0 dx^c \otimes dx^0
\\ = & \partial_c x^{-2}\hat{T}_0 dx^c \otimes dx^0
\\ = & x^{-2} \partial_c\hat{T}_0dx^c \otimes dx^0 
\\ = & (x\partial_c\hat{T}_0) \frac{dx^c}{x}\otimes \frac{dx^0}{x^2},
\\ & \Gamma^d_{bc}T_d = \Gamma^d_{0c}T_d,\\
& \Gamma^d_{0c} = \frac{1}{2} g^{dr}( \partial_cg_{r0}+\partial_xg_{rc}-\partial_rg_{0c} ).
\end{align*}
The argument about the Christoffel symbol is the same as the previous case after interchange the indices $c$ and $b$ (0 here), and we have $\Gamma^d_{0c}T_d dx^c \otimes dx^0 = O(x)\frac{dx^c}{x}\otimes\frac{dx}{x^2}$.

Case 4: $c \neq 0, b\neq 0$.
\begin{align*}
& \partial_cT_b dx^c\otimes dx^b 
\\ = & \partial_c (x^{-1}\hat{T}_b) dx^c\otimes dx^b
\\ = & (x\partial_c\hat{T}_b )\frac{dx^c}{x}\otimes \frac{dx^b}{x},\\
& \Gamma^d_{bc} = \frac{1}{2} g^{dr}( \partial_cg_{rb}+\partial_bg_{rc}-\partial_rg_{bc} )= O(1).
\end{align*}
Next we explain the last line. Since $g_{ij} = O(x^{-2})$, when $r \neq 0$, terms in the bracket are $O(x^{-2})$ altogether. Since $g^{dr}=O(x^2)$ (including the $O(x^4)$ case), this part gives $O(1)$ contribution. When $r = 0$, we have $g^{dr}=O(x^4)$ by (\ref{dual_metric}), and terms in the bracket is of order at most $O(x^{-3})$, hence this part gives $O(x)$ contribution. Combining two cases $r=0$ and $r\neq 0$, we have $\Gamma^d_{bc} = O(1)$, thus $\Gamma^d_{bc}T_ddx^c\otimes dx^b = \sum_dO(1)x\hat{T}_d \frac{dx^b}{x}\otimes \frac{dx^c}{x}$. Thus those terms are $O(x)$ small compared with main terms.

Combining four cases and symmetrize $\nabla$, we have the decomposition
\begin{align}
{d}^{s} = d_0^s + xA, \, A \in \Psi_{\mathrm{sc}}^{0,0}(X;\leftidx{^{\mathrm{sc}}}{T^*X},\mathrm{Sym}^{2,\mathrm{sc}}T^*X), \label{ds_decomposition}
\end{align}
thus as a semiclassical foliation scattering operator, $d^s$ has the same principal symbol as the exterior differential.

\subsection{$d^s$ as a semiclassical foliation 1-cusp operator}
\label{section_ds_1c}
Next we compute the principal symbol of $d^s$ sending 1-cusp one forms to 1-cusp 2-tensors. 
We consider the contribution introduced when we change bundles by comparing the basis of the scattering cotangent bundle and the basis of the foliation semiclassical 1-cusp cotangent bundle.
Initially ${d}^{s}$ is a first order differential operator sending sections of $\leftidx{^{\mathrm{sc}}}{T^*X}$ to sections of $\mathrm{Sym}^{2,\mathrm{sc}}T^*X$. The standard principal symbol of ${d}^{s}$ is tensoring with the covector at which the principal symbol is evaluated, which coincides in the first order with that when we consider it as a first order differential operator sending sections of $\leftidx{_{h,\mathcal{F}}^{\;\;\mathrm{1c}}}{T^*X}$ to sections of $\mathrm{Sym}_{h,\mathcal{F}}^{2,\mathrm{1c}}T^*X$, which is the symmetric part of $\leftidx{_{h,\mathcal{F}}^{\;\;\mathrm{1c}}}{T^*X} \otimes \leftidx{_{h,\mathcal{F}}^{\;\;\mathrm{1c}}}{T^*X}$. 
Combining with (\ref{ds_decomposition}), the zeroth order part introduced by this bundle change
\emph{forms a matrix with $x$ times smooth coefficients in the local basis of $\mathrm{hom}(^{\mathrm{sc}}T^*X,\mathrm{Sym}^{2,\mathrm{sc}}T^*X)$.} Next we consider their contribution in terms of $\mathrm{hom}(\leftidx{_{h,\mathcal{F}}^{\;\;\mathrm{1c}}}{T^*X},\mathrm{Sym}_{h,\mathcal{F}}^{2,\mathrm{1c}}T^*X)$.
With $(\frac{dx}{x^2})^2 := \frac{dx}{x^2} \otimes \frac{dx}{x^2}$, $\frac{dx}{x^2}\frac{dy_i}{x} := \frac{1}{2}(\frac{dx}{x^2} \otimes\frac{dy_i}{x}+ \frac{dy_i}{x}\otimes \frac{dx}{x^2})$, the local basis of $\mathrm{hom}(\leftidx{^{\mathrm{sc}}}{T^*X},\mathrm{Sym}^{2,\mathrm{sc}}T^*X)$ is:
\begin{align}
\begin{split}
& \frac{dx^2}{x^4} \otimes (x^2\partial_x), \,
\frac{dx^2}{x^4} \otimes (x \partial_{y_j}), \,
\frac{dx}{x^2}\frac{dy_i}{x} \otimes (x^2 \partial_x), \\
& \frac{dx}{x^2} \frac{dy_i}{x} \otimes (x^2\partial_{y_j}),  \,
\frac{dy_kdy_i}{x^2} \otimes (x^2\partial_x), \,
\frac{dy_kdy_i}{x^2} \otimes (x\partial_{y_j}),
\end{split} \label{basis_sc}
\end{align}
and the local basis of $\mathrm{hom}(\leftidx{_{h,\mathcal{F}}^{\;\;\mathrm{1c}}}{T^*X},\mathrm{Sym}_{h,\mathcal{F}}^{2,\mathrm{1c}}T^*X)$ is:
\begin{align}
\begin{split}
& \frac{dx^2}{h^2x^6} \otimes (hx^3\partial_x), \,
\frac{dx^2}{h^2x^6} \otimes (h^{1/2}x^2\partial_{y_j}), \,
\frac{dx}{hx^3}\frac{dy_i}{h^{1/2}x} \otimes ( hx^3 \partial_x), \\
& \frac{dx}{hx^3} \frac{dy_i}{h^{1/2}x} \otimes (h^{1/2}x\partial_{y_j}),  \,
\frac{dy_kdy_i}{hx^2} \otimes (hx^3\partial_x), \,
\frac{dy_kdy_i}{hx^2} \otimes (h^{1/2}x\partial_{y_j}).
\end{split} \label{basis_1c}
\end{align}
Comparing the power of $h$ and $x$, in terms of basis in (\ref{basis_1c}), basis in  (\ref{basis_sc}) are smooth and vanish at $\{x=h=0\}$ to orders in following table
\begin{center}
\begin{tabular}{|c|c|}
$\frac{dx^2}{x^4} \otimes (x^2\partial_x)$ & $hx$\\
$\frac{dx^2}{x^4} \otimes (x \partial_{y_j})$ & $h^{3/2}x$\\ 
$\frac{dx}{x^2}\frac{dy_i}{x} \otimes ( x^2 \partial_x)$ & $ h^{1/2} $\\
$\frac{dx}{x^2} \frac{dy_i}{x} \otimes (x\partial_{y_j})$ & $hx$   \\ 
$\frac{dy_kdy_i}{x^2} \otimes (x^2\partial_x)$ & $x^{-1}$\\
$\frac{dy_kdy_i}{x^2} \otimes (x\partial_{y_j})$ & $h^{1/2}$
\end{tabular}
\end{center}
Taking the overall $x$ factor into consideration, they vanish to order $hx^2,h^{3/2}x,$ $h^{1/2}x,hx^2,O(1),h^{1/2}x$ respectively.
Thus the only non-trivial contribution of $d^s-d_{g_{h,\mathrm{1c}}}^s$ is from $\frac{dy_kdy_i}{x^2} \otimes (x^2\partial_x)$ component.
Again using a computation similar to the one for $d^s$ as a scattering operator, the gradient with respect to $g_{h,\mathrm{1c}}$, as a semiclassical 1-cusp operator, has principal symbol
\begin{align*}
\begin{pmatrix}
\xi_{\mathrm{1c}} && 0\\
\eta _{\mathrm{1c}}  && 0 \\
0 && \xi_{\mathrm{1c}} \\
0 && \eta_{\mathrm{1c}} \otimes 
\end{pmatrix},
\end{align*}
where the matrix is acting on matrices of the form $\begin{pmatrix} \xi_1 \\ \eta_1 \end{pmatrix}$ representing $\xi_1 \frac{dx}{hx^3} + \eta_1 \frac{dy}{h^{1/2}x}$. Recall that the effect of conjugation by $e^{-\frac{\F\Phi}{h}}$ is replacing $\xi_{\mathrm{1c}}$ by $\xi_{\mathrm{1c}}-i\F$, thus after taking the error term above, symmetrization and conjugation into consideration, we know that the principal symbol of $d_h^s=e^{\frac{\F}{2hx^2}}d^se^{-\frac{\F}{2hx^2}}$ viewed as an operator between 1-cusp sections has the form
\begin{align*}
\begin{pmatrix}
\xi_{\mathrm{1c}}-i\F && 0 \\
\frac{1}{2}\eta_{\mathrm{1c}} \otimes   && \frac{1}{2}(\xi_{\mathrm{1c}}-i\F)\\
\frac{1}{2}\eta_{\mathrm{1c}}\otimes    && \frac{1}{2}(\xi_{\mathrm{1c}}-i\F)\\
b_s && \eta_{\mathrm{1c}} \otimes_s 
\end{pmatrix},
\end{align*}
where $b_s$ essentially only plays role at the boundary by acting on $\frac{dx}{hx^3}$ to produce 2-tensors. In particular, $b_s$ has 0 differential order. As a consequence, its adjoint $\delta_h^s$ with respect to the metric given in (\ref{metric_sc_1c_1}), acting on symmetric 2-tensors, has principal symbol 
\begin{align}
\begin{pmatrix}
\xi_{\mathrm{1c}} +i\F && \frac{1}{2}\iota_{\eta_{\mathrm{1c}}}  && \frac{1}{2}\iota_{\eta_{\mathrm{1c}}} &&  \la b_s,\cdot \ra \\
0 && \frac{1}{2}(\xi_{\mathrm{1c}} +i\F) && \frac{1}{2}(\xi_{\mathrm{1c}} +i\F) &&  \iota_{\eta_{\mathrm{1c}}} 
\end{pmatrix}, \label{deltas_symbol}
\end{align} 
where $\iota_{\eta_{\mathrm{1c}}}^s=\frac{1}{2}(\eta_{\mathrm{1c},i} \delta_{lj}+\eta_{\mathrm{1c},j}\delta_{il})$ on the lower right corner is replaced by $\iota_{\eta_{\mathrm{1c}}}$ since we only consider symmetric 2-tensors, on which they have the same action. Summarizing results of two parts of this section we have:
\begin{prop}
On functions, the operator $\Delta_{h,\F,s}:=\delta^s_{h,\F} d^s_{h,\F}\in \mathrm{Diff}_{h,\mathrm{1c}}^{2,0}(X)$ has principal symbol 
\begin{align*}
\begin{pmatrix}
\xi_{\mathrm{1c}}+i\F && \iota_{\eta_{\mathrm{1c}}}
\end{pmatrix}
\begin{pmatrix}
 \xi_{\mathrm{1c}}-i\F \\  \eta_{\mathrm{1c}}
 \end{pmatrix} 
= \xi_{\mathrm{1c}}^2+\F^2+|\eta_{\mathrm{1c}}|^2.
\end{align*} 
On one forms,  $\Delta_{h,\F,s}:= \delta^s_{h,\F} d^s_{h,\F} \in \mathrm{Diff}_{h,\mathrm{1c}}^{2,0}(X, \leftidx{_{h,\mathcal{F}}^{\;\;\mathrm{1c}}}{T^*X},\leftidx{_{h,\mathcal{F}}^{\;\;\mathrm{1c}}}{T^*X})$ has principal symbol 
\begin{align}
\scriptsize
\begin{split}
& \begin{pmatrix}
\xi_{\mathrm{1c}} +i\F && \frac{1}{2}\iota_{\eta_{\mathrm{1c}} }  && \frac{1}{2}\iota_{\eta_{\mathrm{1c}}} &&  \la b_s,\cdot \ra \\
0 && \frac{1}{2}(\xi_{\mathrm{1c}} +i\F) && \frac{1}{2}(\xi_{\mathrm{1c}} +i\F) &&  \iota_{\eta_{\mathrm{1c}}} 
\end{pmatrix}
\begin{pmatrix}
\xi_{\mathrm{1c}}-i\F && 0 \\
\frac{1}{2}\eta_{\mathrm{1c}} \otimes  && \frac{1}{2}(\xi_{\mathrm{1c}}-i\F)\\
\frac{1}{2}\eta_{\mathrm{1c}}\otimes   && \frac{1}{2}(\xi_{\mathrm{1c}}-i\F)\\
b_s && \eta_{\mathrm{1c}} \otimes_s 
\end{pmatrix}
\\ = &
\begin{pmatrix}
\xi_{\mathrm{1c}}^2+\F^2+\frac{1}{2}\eta_{\mathrm{1c}}^2 && \frac{1}{2}(\xi_{\mathrm{1c}}-i\F)\iota_{\eta_{\mathrm{1c}}} \\
\frac{1}{2}(\xi_{\mathrm{1c}}+i\F)\eta_{\mathrm{1c}} \otimes && \frac{1}{2}(\xi_{\mathrm{1c}}^2+\F^2)+ \iota_{\eta_{\mathrm{1c}}} \eta_{\mathrm{1c}}\otimes_s
\end{pmatrix}
 + \begin{pmatrix}
\la b_s,\cdot \ra b_s && \la b_s,\cdot \ra \eta_{\mathrm{1c}}\otimes_s \\
\iota_{\eta_{\mathrm{1c}}}b_s && 0
\end{pmatrix}.
\end{split}
\label{Delta_symbol}
\end{align}
\label{prop_delta_symbol}
\end{prop}

\subsection{$d^s$ as a semiclassical foliation scattering operator near $\Sigma_{x_0}$}
\label{section_ds_sc}
In this section we consider $d^s$ and the modified Laplacian $\Delta_{h,\F,s}$ near the artificial boundary as semiclassical foliation scattering operators and compute their symbols. 
This is similar to the argument in Section \ref{section_ds_1c}: we consider the contribution introduced when we change bundles by comparing the basis of the scattering cotangent bundle and the basis of the foliation semiclassical 1-cusp cotangent bundle. 
As a first order differential operator sending sections of
$\leftidx{^{\mathrm{sc}}}{T^*X}$ to sections of
$\mathrm{Sym}^{2,\mathrm{sc}}T^*X$, the standard principal symbol of
${d}^{s}$ is tensoring with the covector at which the principal symbol
is evaluated, which coincides with that when we consider it as a first order differential operator sending sections of $\leftidx{_{h,\mathcal{F}}^{\;\;\mathrm{sc}}}{T^*X}$ to sections of $\mathrm{Sym}_{h,\mathcal{F}}^{2,\mathrm{sc}}T^*X$, which is the symmetric part of $\leftidx{_{h,\mathcal{F}}^{\;\;\mathrm{sc}}}{T^*X} \otimes \leftidx{_{h,\mathcal{F}}^{\;\;\mathrm{sc}}}{T^*X}$. 

We emphasize that the scattering structure of
$\leftidx{^{\mathrm{sc}}}{T^*X}$ in Section~\ref{section_Chris} refers
to using $\Sigma_0$ as the boundary while here
$\leftidx{_{h,\mathcal{F}}^{\;\;\mathrm{sc}}}{T^*X}$  refers to using
$\Sigma_{x_0}$ as the boundary. The operator  $d^s$ here is the symmetric differential with respect to a smooth metric.
The zeroth order part introduced by the bundle change from the
standard smooth bundles to the new (local, near $\Sigma_{x_0}$) scattering bundle
\emph{forms a matrix with smooth coefficients in the local basis of $\mathrm{hom}(^{\mathrm{sc}}T^*X,\mathrm{Sym}^{2,\mathrm{sc}}T^*X)$}, since now there is no gain for the zeroth order part as in (\ref{ds_decomposition}) in terms of $(x_0-x)$.
Next we consider their contribution in terms of $\mathrm{hom}(\leftidx{_{h,\mathcal{F}}^{\;\;\mathrm{sc}}}{T^*X},\mathrm{Sym}_{h,\mathcal{F}}^{2,\mathrm{sc}}T^*X)$.
Recall that the basis of
$\mathrm{hom}(^{\mathrm{sc}}T^*X,\mathrm{Sym}^{2,\mathrm{sc}}T^*X)$ is
given by (\ref{basis_sc}) with $\rho=x_0-x$ in place of $x$, while
the local basis (near $\Sigma_{x_0}$) of
$\mathrm{hom}(_{h,\mathcal{F}}^{\;\;\mathrm{sc}}T^*X,\mathrm{Sym}_{h,\mathcal{F}}^{2,\mathrm{sc}}T^*X)$,
using again $\rho=x_0-x$ for brevity, is:
\begin{align}
\begin{split}
& \frac{d\rho^2}{h^2\rho^4} \otimes (h\rho^2\partial_\rho), \,
\frac{d\rho^2}{h^2\rho^4} \otimes (h^{1/2} \rho \partial_{y_j}), \,
\frac{d\rho}{h\rho^2}\frac{dy_i}{h^{1/2}\rho} \otimes (h\rho^2 \partial_\rho), \\
& \frac{d\rho}{h\rho^2} \frac{dy_i}{h^{1/2}\rho} \otimes (h \rho^2\partial_{y_j}),  \,
\frac{dy_kdy_i}{h\rho^2} \otimes (h\rho^2\partial_\rho), \,
\frac{dy_kdy_i}{h\rho^2} \otimes (h^{1/2}\rho\partial_{y_j}).
\end{split} \label{basis_sc_h}
\end{align}
Comparing the power of $h$ (powers of $\rho$ are all the same in this case), in terms of basis in (\ref{basis_sc_h}), basis in  (\ref{basis_sc}) are smooth and vanish at $\{\rho=h=0\}$ to orders in following table
\begin{center}
\begin{tabular}{|c|c|}
$\frac{d\rho^2}{\rho^4} \otimes (\rho^2\partial_\rho)$ & $h$\\
$\frac{d\rho^2}{\rho^4} \otimes (\rho \partial_{y_j})$ & $h^{3/2}$\\ 
$\frac{d\rho}{\rho^2}\frac{dy_i}{\rho} \otimes ( \rho^2 \partial_\rho)$ & $ h^{1/2}$\\
$\frac{d\rho}{\rho^2} \frac{dy_i}{\rho} \otimes (\rho\partial_{y_j})$ & $h$   \\ 
$\frac{dy_kdy_i}{\rho^2} \otimes (\rho^2\partial_\rho)$ & $O(1)$\\
$\frac{dy_kdy_i}{\rho^2} \otimes (\rho\partial_{y_j})$ & $h^{1/2}$
\end{tabular}
\end{center}
Again the only non-trivial contribution is from $\frac{dy_kdy_i}{h\rho^2} \otimes (h\rho^2\partial_\rho)$. 

Recall that $\Phi=\frac{1}{x_0-x}$ near $\Sigma_{x_0}$, the effect of the conjugation near $\Sigma_{x_0}$, is derived from
\begin{align*}
e^{-\frac{\F}{h(x_0-x)}} h(x_0-x)^2D_x e^{\frac{\F}{h(x_0-x)}} = h(x_0-x)^2D_x -i\F,
\end{align*}
which tells us the conjugation is effectively replacing $\xi_{\mathrm{sc}}$ by $\xi_{\mathrm{sc}}-i\F$.
Thus $d^s$ when viewed as a semiclassical foliation scattering operator near $\Sigma_{x_0}$ has principal symbol
\begin{align}
\begin{pmatrix}
\xi_{\mathrm{sc}}-i\F && 0 \\
\frac{1}{2}\eta_{\mathrm{sc}} \otimes   && \frac{1}{2}(\xi_{\mathrm{sc}}-i\F)\\
\frac{1}{2}\eta_{\mathrm{sc}}\otimes    && \frac{1}{2}(\xi_{\mathrm{sc}}-i\F)\\
\tilde{b}_s && \eta_{\mathrm{sc}} \otimes_s 
\end{pmatrix},
\label{ds_symbol_sc}
\end{align}
where $\tilde{b}_s$ only plays role at $\Sigma_{x_0}$ by acting on $\frac{dx}{h(x_0-x)^2}$ to produce 2-tensors. In particular, $\tilde{b}_s$ has 0 differential order. Its adjoint $\delta_h^s$ with respect to the metric in (\ref{metric_sc_1c_2}), acting on symmetric 2-tensors, has principal symbol 
\begin{align}
\begin{pmatrix}
\xi_{\mathrm{sc}} +i\F && \frac{1}{2}\iota_{\eta_{\mathrm{sc}} }  && \frac{1}{2}\iota_{\eta_{\mathrm{sc}}} &&  \la \tilde{b}_s,\cdot \ra \\
0 && \frac{1}{2}(\xi_{\mathrm{sc}} +i\F) && \frac{1}{2}(\xi_{\mathrm{sc}} +i\F) &&  \iota_{\eta_{\mathrm{sc}} } 
\end{pmatrix}, \label{deltas_symbol_sc}
\end{align} 
where $\iota_{\eta_{\mathrm{sc}}}^s=\frac{1}{2}(\eta_{\mathrm{sc},i} \delta_{lj}+\eta_{\mathrm{sc},j}\delta_{il})$ on the lower right corner is replaced by $\iota_{\eta_{\mathrm{sc}}}$ since we are acting on symmetric 2-tensors. Taking product of (\ref{ds_symbol_sc}), (\ref{deltas_symbol_sc}) and summarizing, we have shown:
\begin{prop}
On functions, the operator $\Delta_{h,\F,s}:=\delta^s_{h,\F} d^s_{h,\F}\in \mathrm{Diff}_{h,\mathrm{sc}}^{2,0}(X)$ has principal symbol 
\begin{align*}
\begin{pmatrix}
\xi_{\mathrm{sc}}+i\F && \iota_{\eta_{\mathrm{sc}}}
\end{pmatrix}
\begin{pmatrix}
 \xi_{\mathrm{sc}}-i\F \\  \eta_{\mathrm{sc}}
 \end{pmatrix} 
= \xi_{\mathrm{sc}}^2+\F^2+|\eta_{\mathrm{sc}}|^2.
\end{align*} 
On one forms,  $\Delta_{h,\F,s}:= \delta^s_{h,\F} d^s_{h,\F} \in \mathrm{Diff}_{h,\mathrm{sc}}^{2,0}(X; \leftidx{_{h,\mathcal{F}}^{\;\;\mathrm{sc}}}{T^*X},\leftidx{_{h,\mathcal{F}}^{\;\;\mathrm{sc}}}{T^*X})$ has principal symbol 
\begin{align}
\scriptsize
\begin{split}
& \begin{pmatrix}
\xi_{\mathrm{sc}} +i\F && \frac{1}{2}\iota_{\eta_{\mathrm{sc}} }  && \frac{1}{2}\iota_{\eta_{\mathrm{sc}}} &&  \la \tilde{b}_s,\cdot \ra \\
0 && \frac{1}{2}(\xi_{\mathrm{sc}} +i\F) && \frac{1}{2}(\xi_{\mathrm{sc}} +i\F) &&  \iota_{\eta_{\mathrm{sc}}} 
\end{pmatrix}
\begin{pmatrix}
\xi_{\mathrm{sc}}-i\F && 0 \\
\frac{1}{2}\eta_{\mathrm{sc}} \otimes  && \frac{1}{2}(\xi_{\mathrm{sc}}-i\F)\\
\frac{1}{2}\eta_{\mathrm{sc}}\otimes   && \frac{1}{2}(\xi_{\mathrm{sc}}-i\F)\\
\tilde{b}_s && \eta_{\mathrm{sc}} \otimes_s 
\end{pmatrix}
\\ = &
\begin{pmatrix}
\xi_{\mathrm{sc}}^2+\F^2+\frac{1}{2}\eta_{\mathrm{sc}}^2 && \frac{1}{2}(\xi_{\mathrm{sc}}-i\F)\iota_{\eta_{\mathrm{sc}}} \\
\frac{1}{2}(\xi_{\mathrm{sc}}+i\F)\eta_{\mathrm{sc}} \otimes && \frac{1}{2}(\xi_{\mathrm{sc}}^2+\F^2)+ \iota_{\eta_{\mathrm{sc}}} \eta_{\mathrm{sc}}\otimes_s
\end{pmatrix}
 + \begin{pmatrix}
\la \tilde{b}_s,\cdot \ra \tilde{b}_s && \la \tilde{b}_s,\cdot \ra \eta_{\mathrm{sc}}\otimes_s \\
\iota_{\eta_{\mathrm{sc}}}\tilde{b}_s && 0
\end{pmatrix}.
\end{split}
\label{Delta_symbol_sc}
\end{align}
\label{prop_delta_symbol_sc}
\end{prop}

\section{The modified normal operator}\label{sec:normal}
In this section we consider the membership and ellipticity of the modified normal operator restricted to the kernel of $\delta_{h,\F}^s$.
Recall that $I$ is the X-ray transform defined by (\ref{xray_defn}) and $\Phi$ is defined in (\ref{Phi_defn}), we define the modified normal operator $N_{h,\F}$ of $I$ by
\begin{align*}
N_{h,\F} = e^{-\frac{\F\Phi}{h}}L \tilde{\chi} I e^{\frac{\F \Phi}{h}},
\end{align*}
where $L$ the `adjoint' of $I$ defined for one forms by
\begin{align*}
(Lw)(z) =  x^2\int_{S_zM} w(\gamma_{x,y,\lambda,\omega})g_{\mathrm{sc},\mathrm{1c},h}(\dot{\gamma}_{x,y,\lambda,\omega}(0))d\lambda d\omega,
\end{align*}
and for 2-tensors by
\begin{align*}
(Lw)(z) = h x^4\int_{S_zM} w(\gamma_{x,y,\lambda,\omega})g_{\mathrm{sc},\mathrm{1c},h}(\dot{\gamma}_{x,y,\lambda,\omega}(0)) \otimes g_{\mathrm{sc},\mathrm{1c},h}(\dot{\gamma}_{x,y,\lambda,\omega}(0)) d\lambda d\omega,
\end{align*}
where $g_{\mathrm{sc},\mathrm{1c},h}$ is the metric defined using
(\ref{metric_sc_1c_1}) and (\ref{metric_sc_1c_2}). As we mentioned,
$I$ sends functions on $M$ to functions on the collections of
geodesics on $M$, which is identified as $TM$ by identifying the
starting point and the tangent vector of a geodesic with itself. Here
the $x^2$ and $hx^4$ pre-factors are introduced to cancel factors
introduced by
$g_{\mathrm{sc},\mathrm{1c},h}(\dot{\gamma}_{x,y,\lambda,\omega}(0))$
and its tensor powers, similarly to how powers of $x$ were used in \cite{stefanov2018inverting}. In fact, it would be conceptually more clear to use $hx^2$ and $h^2x^4$ respectively, but our computation, especially the part after we use the rescaled variables in the symbol computation, shows that results are more compact if we have 1 order less $h$ power in our definitions.

Next we discuss quantities in the definition of $L$ in more detail. Recall (\ref{conic_metric_2}), the asymptotic conic metric is
\begin{align*}
g = x^{-4}dx^2 + x^{-2}\tilde g,
\end{align*}
with $\tilde g|_{x=0}$ on the cross section being the asymptotic link
metric $g_0$.
The computation in  \cite{Melrose1994} shows that, the Hamilton vector filed associated to the dual metric function $G$ of the asymptotic conic metric $g$ is 
\begin{align*}
H_G = 2x( (\xi_{\mathrm{sc}} x\partial_x+\eta_{\mathrm{sc}}\cdot \partial_{\eta_{\mathrm{sc}}}) - |\eta_{\mathrm{sc}}|^2\partial_{\eta_{\mathrm{sc}}}+\frac{1}{2}H_{G_1}  +xV),
\end{align*}
where $G_1$ is the dual metric function of $g_1$ and $V$ is a vector field tangent to the boundary $\{x=0\}$.
The geodesic is given by (page 28 of \cite{zachos2022inverting})
\begin{align}
\gamma_{x,y,\lambda,\omega}(t) = (x+x(\lambda t + \alpha t^2+t^3\Gamma^{(1)}),y+\omega t+t^2 \Gamma^{(2)}),
\label{geodesic}
\end{align}
where $\Gamma^{(1)},\Gamma^{(2)}$ are smooth functions in $x,y,\lambda,\omega,t$.
Taking derivative with respect to $t$, the tangent vector at $x'$ is
\begin{align*}
\lambda'x\partial_x+\omega'\partial_y=\dot{\gamma}_{x,y,\lambda,\omega}(t) = (x\lambda + 2 x \alpha t+O(xt^2))\partial_x+(\omega+O(t))\partial_y.
\end{align*}
Setting $t=0$, the tangent vector at $x$ is
\begin{align}
x\lambda \partial_x + \omega\partial_y.  \label{tangent}
\end{align}
A key point in \cite[Section~3.1]{zachos2022inverting} as well as in the earlier
works is the negativity of $\alpha$ at $\lambda=0$; this precisely
corresponds to a strict (definite) concavity statement on the level sets of $x$ from the
sublevel sets.

The cutoff $\tilde{\chi}$ and the weight function $\Phi$ are also
designed to depend on the geometric/analytic setting. We  choose
$\Phi(x)=-\frac{1}{2x^2}$ in the analytic 1-cusp setting (which is
geometrically scattering) and $\Phi(x)=\frac{1}{x}$ in the scattering
setting (which is geometrically smooth). Choose $\tilde{\chi}$ to be
\begin{align*}
\tilde{\chi} = \chi(x^{1/2} \sqrt{\Phi'} (h^{1/2}|\alpha|^{1/2} )^{-1} ),
\end{align*}
where $\chi \in C_c^\infty(\R)$ is non-negative and identically 1 near 0 and $\alpha$ is introduced in (\ref{geodesic}). 

Membership of $N_{h,\F}$ in ${\Psi}_{\mathrm{sc},h,\mathcal{F}}^{-1,0}(X;\leftidx{_{h,\mathcal{F}}^{\;\;\mathrm{sc}}}{T^*X},\leftidx{_{h,\mathcal{F}}^{\;\;\mathrm{sc}}}{T^*X})$ for the one form case and ${\Psi}_{\mathrm{sc},h,\mathcal{F}}^{-1,2}(X;\mathrm{Sym}_{h,\mathcal{F}}^{2,\mathrm{sc}}T^*X,\mathrm{Sym}_{h,\mathcal{F}}^{2,\mathrm{sc}}T^*X)$ for the two tensor case when localized near the artificial boundary $\{x=x_0\}$ follows from Proposition 3.3 of \cite{vasy2020semiclassical} and the argument in Proposition 3.1 of \cite{stefanov2018inverting}, which generalizes Proposition 3.3 of \cite{uhlmann2016inverse} to one form and two tensor cases, just notice that $x$ in \cite{stefanov2018inverting} is $(x_0-x)$ now, and the powers of it in (2.1)(2.2) of \cite{stefanov2018inverting} are encoded in the second index of ${\Psi}_{\mathrm{sc},h,\mathcal{F}}^{*,*}$. More detailed dicussion is given in the proof of Proposition \ref{prop_sc}.

The membership of $N_{h,\F}$ localized near the conic infinity $x=0$ in ${\Psi}_{\mathrm{1c},h,\mathcal{F}}^{-1,-1}$ follows from Theorem 3.1 of \cite{zachos2022inverting}. The only difference between our operator and operator therein is the additional tensorial factor, which is a smooth endomorphism between scattering and 1-cusp bundles when localize to each ends respectively and does not affect the pseudodifferential property. The power of $x$ and $h$ introduced by those tensorial factors are cancelled by those in the definition of $L$. So our proof focuses on the ellipticity.
\begin{prop}
Let $\F>0$ for one forms, and $\F$ is sufficiently large for two tensors and $\Omega_{x_0} = \{ x \leq x_0 \}$ with $x_0$ small.
$N_{h,\F}$ is an operator in
${\Psi}_{\mathrm{1c},h,\mathcal{F}}^{-1,0}$ and it is elliptic on the
kernel of $\delta_{h,\F}$ in the semiclassical foliation 1-cusp algebra.
 \label{prop_conic}
\end{prop}
Because the statement involves ellipticity in various senses for both
one forms and 2-tensors (differential, semiclassical, boundary symbol), we decompose the proof into several lemmas, whose combination proves Proposition \ref{prop_conic}. 

\subsection{Ellipticity on one forms}
In this part we show ellipticity of $N_{h,\F}$ on one forms in both
1-cusp differential and semiclassical and boundary sense. We first
start with the differential behavior, corresponding to `infinite
points' in the semiclassical foliation 1-cusp cotangent bundle (fiber
infinity for the compactification), and then we turn to finite points
(semiclassical and boundary behavior) as the arguments are somewhat
different in the two cases.
\begin{lmn}
For $\F>0$,
$N_{h,\F}$ is an operator in
${\Psi}_{\mathrm{1c},h,\mathcal{F}}^{-1,-1}$ acting on one forms and it
is elliptic on the kernel of $\delta_{h,\F}$ in the 1-cusp algebra
differential sense near $\Sigma_0$.
 \label{lmn_oneform_infinity}
\end{lmn}

\begin{proof}
First we compute the tensorial factor $g_{h,\mathrm{1c}}(\dot{\gamma},\cdot)\otimes \iota_{\dot{\gamma}}$
in the Schwartz kernel of $N_{h,\F}$ at the front face $\{x=0\}$. Introduce the rescaled variables:
\begin{align*}
\hat{\lambda} = \frac{\lambda}{h^{1/2}x}, \hat{t} = \frac{t}{h^{1/2}x}.
\end{align*}
The semiclassical 1-cusp metric applied to (\ref{tangent}) is:
\begin{align*}
g_{h,\mathrm{1c}}(x\lambda \partial_x + \omega\partial_y) & = \lambda \frac{dx}{h^2x^5}+\frac{g_1(\omega\partial_y)}{hx^2}\\
& =  h^{-1/2}x^{-1}(\hat{\lambda}\frac{dx}{hx^3} + \frac{g_1(\omega\partial_y)}{h^{1/2}x}).
\end{align*}
Similarly, the second factor is
\begin{align*}
\iota_{\dot{\gamma}} = h^{-1/2}x^{-1}( (\hat{\lambda}+2\alpha\hat{t})hx^3\partial_x + \omega h^{1/2}x\partial_y).
\end{align*}
So the tensorial factor near $\Sigma_0$ is,
\begin{align*}
\tilde{E}_c = h^{-1}(\hat{\lambda}\frac{dx}{hx^3} + \frac{g_1(\omega\partial_y)}{h^{1/2}x}) \otimes ( (\hat{\lambda}+2\alpha\hat{t})hx^3\partial_x + \omega h^{1/2}x\partial_y),
\end{align*}
where we have incorporated the $x^2$ factor in the definition of $L$ on one forms. 
The modified normal operator $N_{h,\F}$ acts on functions by
\begin{align*}
N_{h,\F}u(x,y) =  (2\pi)^{-n} h^{-n-1/2} & \int e^{ix^{-3}\tilde{\xi}_{\mathrm{1c}}\frac{(\gamma^{(1)}_{x,y,\lambda,\omega}(t)-x)}{h}}e^{ix^{-1}\tilde{\eta}_{\mathrm{1c}}\cdot \frac{(\gamma^{(2)}_{x,y,\lambda,\omega}(t)-y)}{h^{1/2}}} 
 \\
& a_h(x,y,\tilde{\xi}_{\mathrm{1c}} , \tilde{\eta}_{\mathrm{1c}}) (x')^{-n-2} dx'dy' d \tilde{\xi}_{\mathrm{1c}} d\tilde{\eta}_{\mathrm{1c}},
\end{align*}
where $a_h$ is its standard left symbol.
The Schwartz kernel of $N_{h,\F}$ relative to the density $|dz'|=|dx'dy'|$ is 
\begin{align}
\begin{split}
K_{N_{h,\F}}(x,y,x',y')= &\int e^{-\F\Phi(x)/h}e^{\F\Phi(x(\gamma_{x,y,\lambda,\omega}(t)))/h}\tilde{\chi}(x,y,\lambda/xh^{1/2},\omega)
\\& \delta(z'-\gamma_{z,\lambda,\omega}(t))\tilde{E}_c (\gamma^{(1)_{z,\lambda,\omega}(t)})^{-1} dtd\lambda d\omega\\
= & (2\pi)^{-n}h^{-n/2-1/2} \int e^{-\F\Phi(x)/h}e^{\F\Phi(x(\gamma_{x,y,\lambda,\omega}(t)))/h}
\\&\tilde{\chi}(x,y,\lambda/xh^{1/2},\omega)\tilde{E}_c e^{-i\xi'\cdot\frac{(x'-\gamma^{(1)}_{x,y,\lambda,\omega}(t))}{h}-i\eta'\cdot \frac{(y'-\gamma^{(2)}_{x,y,\lambda,\omega}(t))}{h^{1/2}}} 
\\&(\gamma^{(1)_{z,\lambda,\omega}(t)})^{-1}dtd\lambda d\omega d\zeta'.
\end{split}
\label{kernel1}
\end{align}
The principal symbol of $N_{h,\F}$ is the leading part of the semiclassical inverse Fourier transform of $K_{N_{h,\F}}$ in $z'$ evaluated at $\zeta$, which gives a $(2\pi)^n\delta(\zeta-\zeta')$ factor, effectively replacing $\zeta'$ by $\zeta$ in the integrand and annihilate the $\zeta'-$integration. We obtain
\begin{align}
\begin{split}
a_h(x,y,\lambda,\omega)= & \int e^{-\F\Phi(x)/h}e^{\F\Phi(x(\gamma_{x,y,\lambda,\omega}(t)))/h}\tilde{\chi}(x,y,\lambda/xh^{1/2},\omega)\\
&\tilde{E}_c e^{i\xi\cdot\frac{(\gamma^{(1)}_{x,y,\lambda,\omega}(t)-x)}{h}}e^{i\eta\cdot \frac{(\gamma^{(2)}_{x,y,\lambda,\omega}(t)-y)}{h^{1/2}}} dtd\lambda d\omega,
\end{split}
\label{aF}
\end{align}
where $\tilde{\chi}$ has compact support with respect to $\frac{\lambda}{h^{1/2}}$, thus the $\hat{\lambda}-$integral is happening over a compact interval. 

Notice that by the definition of the 1-cusp algebra, our Schwartz kernel is relative to the density $\frac{dx'dy'}{(x')^{n+2}}$, which means it is $(2\pi)^{-n}x^{n+2}$ times the semiclassical foliation Fourier transform in $(x^{-3}\tilde{\xi}_{\mathrm{1c}},x^{-1}\tilde{\eta}_\mathrm{1c})$ of the symbol 
\begin{align*}
(x,y,\tilde{\xi}_{\mathrm{1c}},\tilde{\eta}_{\mathrm{1c}}) \rightarrow
e^{ i(x^{-2}\xi_{\mathrm{1c}}h^{-1} + x^{-1}y\cdot \eta_{\mathrm{1c}}h^{-1/2})}a(x,y,\tilde{\xi}_{\mathrm{1c}},\tilde{\eta}_{\mathrm{1c}}),
\end{align*}
where the $x^{n+2}$ factor comes from scaling the fiber variable to $(x^{-3}\tilde{\xi}_{\mathrm{1c}},x^{-1}\tilde{\eta}_\mathrm{1c})$. Then we invert this Fourier transform and evaluate at $(x^{-3}\tilde{\xi}_{\mathrm{1c}},x^{-1}\tilde{\eta}_\mathrm{1c})$ to obtain
\begin{align*}
a_h(x,y,\xi_{\mathrm{1c}},\eta_{\mathrm{1c}}) = & (2\pi)^nx^{-n-2}e^{ i(-x^{-2}\xi_{\mathrm{1c}}h^{-1} - x^{-1}y\cdot \eta_{\mathrm{1c}}h^{-1/2})}
\\&(\mathcal{F}^{-1}_{h})_{(x',y') \rightarrow (x^{-3}\tilde{\xi}_{\mathrm{1c}},x^{-1}\tilde{\eta}_\mathrm{1c})}K_{N_{h,\F}}(x,y,x',y'),
\end{align*}
where the $\mathcal{F}$ representing (inverse) Fourier transform and the foliation should not cause confusion. 
\begin{align}
\begin{split}
a_h(x,y,\lambda,\omega)= & x^{-n-2}e^{-ix^{-3}\xi_{\mathrm{1c}}h^{-1} -ix^{-1}\eta_{\mathrm{1c}}h^{-1/2} }
 \int e^{-\F\Phi(x)/h}e^{\F\Phi(x(\gamma_{x,y,\lambda,\omega}(t)))/h}
\\& \tilde{\chi}(x,y,\lambda/xh^{1/2},\omega)\tilde{E}_c e^{ix^{-3}\xi_{\mathrm{1c}} \gamma_{z,\lambda,\omega}^{(1)}(t)h^{-1} +ix^{-1}\eta_{\mathrm{1c}}\cdot(\gamma_{z,\lambda,\omega}^{(2)}(t) -y )h^{-1/2} }
\\&(\gamma_{z,\lambda,\omega}^{(1)}(t))^{n+1}dt|d\sigma|\\
&= \int e^{-\F\Phi(x)/h}e^{\F\Phi(x(\gamma_{x,y,\lambda,\omega}(t)))/h}\tilde{\chi}(x,y,\lambda/xh^{1/2},\omega)\tilde{E}_c 
\\&e^{ix^{-3}\xi_{\mathrm{1c}}\frac{(\gamma^{(1)}_{x,y,\lambda,\omega}(t)-x)}{h}}e^{ix^{-1}\eta_{\mathrm{1c}}\cdot \frac{(\gamma^{(2)}_{x,y,\lambda,\omega}(t)-y)}{h^{1/2}}} 
x^{-n-2}(\gamma_{z,\lambda,\omega}^{(1)}(t))^{n+1}dt|d\sigma|.
\end{split}
\end{align}
The imaginary phase is 
\begin{align}
\begin{split}
 & \xi \cdot \frac{(\gamma^{(1)}_{x,y,\lambda,\omega}(t)-x)}{h} + \eta\cdot \frac{(\gamma^{(2)}_{x,y,\lambda,\omega}(t)-y)}{h^{1/2}}\\
=& x^{-3}\xi_{\mathrm{1c}} \frac{(\gamma_{x,y,\lambda,\omega}^{(1)}(t) - x) }{h} + x^{-1}\eta_{\mathrm{1c}}( \frac{(\gamma^{(2)}_{x,y,\lambda,\omega}(t)-y)}{h^{1/2}}  )\\
= & \xi_{\mathrm{1c}}(\hat{\lambda}\hat{t}+\alpha \hat{t}^2+h^{1/2}\hat{t}^2\Gamma^{(1)}(x,y,h^{1/2}\hat{\lambda},\omega,h^{1/2}\hat{t}))
\\& + \eta_{\mathrm{1c}} \cdot ( \omega \hat{t} + h^{1/2}\hat{t}^2\Gamma^{(2)}(x,y,h^{1/2}\hat{\lambda},\omega,h^{1/2}\hat{t}) ).
\end{split}
\label{phase_conic_1}
\end{align}
The damping factor is 
\begin{align*}
-\frac{\F\Phi(x)}{h}+\frac{\F\Phi(x(\gamma_{x,y,\lambda,\omega}(t) ))}{h},
\end{align*}
where $\Phi(x)=-\frac{1}{2x^2}$, thus this damping factor is
\begin{align}
\begin{split}
  & \frac{\F}{2x^2h}-\frac{\F}{2(\gamma^{(1)}_{x,y,\lambda,\omega}(t))^2h}\\
 = & \frac{\F}{2h}( (\gamma^{(1)}_{x,y,\lambda,\omega}(t))^2-x^2)x^{-2}(\gamma^{(1)}_{x,y,\lambda,\omega}(t))^{-2}\\
 = & \frac{\F}{2hx^2}( (1+\lambda t + \alpha t^2 +t^3\Gamma^{(1)})^2-1 )(1+\lambda t + \alpha t^2 +t^3\Gamma^{(1)})^{-2}\\
 = & \frac{\F}{2hx^2}(2+\lambda t+\alpha t^2 +t^3\Gamma^{(1)})(\lambda t+\alpha t^2 +t^3\Gamma^{(1)}) (1+\lambda t + \alpha t^2 +t^3\Gamma^{(1)})^{-2}\\
 = & \F(\hat{\lambda}\hat{t}+\alpha(x,y,xh^{1/2}\hat{\lambda},\omega)\hat{t}^2+\hat{t}^3xh^{1/2}\hat{\Gamma}^{(1)}(x,y,xh^{1/2}\hat{\lambda},\omega,xh^{1/2}\hat{t})),
\end{split}
\label{damp}
\end{align}
where $\hat{\Gamma}^{(1)}$ is smooth with respect to its variables. This means the integrand is Gaussian with respect to $\hat{t}$, hence integrating against it does not affect decay or smoothness properties. Notice that when we rewrite (\ref{aF}) in terms of $\hat{\lambda},\hat{t}$, we obtain an $xh$ factor because of the change of variables, and we combine $h$ with $\tilde{E}_c$ to define
\begin{align*}
E_c=h\tilde{E}_c = (\hat{\lambda}\frac{dx}{hx^3} + \frac{g_1(\omega\partial_y)}{h^{1/2}x}) \otimes ( (\hat{\lambda}+2\alpha\hat{t})hx^3\partial_x + \omega h^{1/2}x\partial_y),
\end{align*}
which is smooth down to $x=0$ and $h=0$.
\begin{align}
\begin{split}
a_h(x,y,\lambda,\omega)= & x\int e^{\F( \hat{\lambda}\hat{t} +\alpha\hat{t}^2+\hat{t}^3h^{1/2}{\Gamma}^{(1)}(x,y,h^{1/2}\hat{\lambda},\omega,h^{1/2}\hat{t}) )} \tilde{\chi}(x,y,\hat{\lambda},\omega) E_c \\
& e^{i(\xi_{\mathrm{1c}}(\hat{\lambda}\hat{t}+\alpha \hat{t}^2+h^{1/2}\hat{t}^2\Gamma^{(1)}) + \eta_{\mathrm{1c}} \cdot ( \omega \hat{t} + h^{1/2}\hat{t}^2\Gamma^{(2)}))} d\hat{t}d\hat{\lambda} d\omega ,
\end{split}
\label{symbol_oneform} 
\end{align}
where we abbreviated the variables $(x,y,h^{1/2}\hat{\lambda},\omega,h^{1/2}\hat{t})$ of $\Gamma^{(1)},\Gamma^{(2)}$. Next we compute $\lambda',\omega'$, i.e. computing $\dot{\gamma}_{x,y,\lambda,\omega}(t)$. Recall the expression of geodesics, we have
\begin{align*}
(\lambda',\omega')=\dot{\gamma}_{x,y,\lambda,\omega}(t) = (\lambda + 2\alpha(x,y,\lambda,\omega)t+O(t^2),\omega+O(t)).
\end{align*}
Recall that $\lambda=h^{1/2}x\hat{\lambda},t=h^{1/2}x\hat{t}$, thus $\lambda'=O(h^{1/2}),\lambda=O(h^{1/2})$, and at the front face $E_c$ is 
\begin{align}
E_c = (\hat{\lambda}\frac{dx}{hx^3} + \frac{g_1(\omega\partial_y)}{h^{1/2}x}) \otimes( (\hat{\lambda}+2\alpha\hat{t})hx^3\partial_x + \omega h^{1/2}x\partial_y).
\label{Ec_1form}
\end{align}
We compute (\ref{symbol_oneform}) by applying the stationary phase lemma. We choose $\tilde{\chi}$ having compact support with respect to its third variable. Thus the integral in $\hat{\lambda},\omega$ is over a compact region, satisfying the condition for the stationary phase lemma. 
We divide the integral in $\hat{t}$ into two parts with one over $\{|\hat{t}| \leq 2\}$ and the other over $\{|\hat{t}| \geq 1\}$ by introducing a partition of unity
 \begin{align*}
\chi_1+\chi_2 =1,
\end{align*}
where $\chi_i \in C^\infty(\R), \, \supp \chi_1 \subset [-2,2], \, \supp \chi_2 \subset (-\infty,-1]\cup[1,\infty)$. In particular, $\chi_1(0)=1$. In $\{ |\hat{t}| \leq 2 \}$ we apply the standard parameter dependent stationary phase lemma. In $\{ |\hat{t}| \geq 1 \}$, the phase in non-stationary and integration by parts shows its Schwartz property.

We first consider the region   $\{|\hat{t}| < 2\}$ by applying the standard parameter dependent stationary phase lemma. In order to compute the semiclassical principal symbol, we may set $h^{1/2}=0$ and the phase becomes
\begin{align*}
\xi_{\mathrm{1c}}(\hat{\lambda}\hat{t}+\alpha(x,y,0,\omega)\hat{t}^2)+\eta_{\mathrm{1c}}\cdot\omega\hat{t}.
\end{align*}
The contribution to the symbol of the part $\{|\hat{t}| < 2\}$ when we take $h^{1/2}=0$ is 
\begin{align*}
a_h(x,y,\xi_{\mathrm{1c}},\eta_{\mathrm{1c}}) = & x \int e^{i(\xi_{\mathrm{1c}}(\hat{\lambda}\hat{t}+\alpha(x,y,0,\omega)\hat{t}^2) + \eta_{\mathrm{1c}}\cdot \omega \hat{t} )} e^{\F(\hat{\lambda}\hat{t} +\alpha(x,y,0,\omega)\hat{t}^2)} \tilde{\chi}(x,y,\hat{\lambda},\omega)
\\& \chi_1(\hat{t}) E_c d\hat{t}d\hat{\lambda}d\omega.
\end{align*}
Error terms caused by taking $h^{1/2}=0$ are of at most $O(xh^{1/2} \la | (\xi_{\mathrm{1c}},\eta_{\mathrm{1c}}) | \ra^{-1} )$ order relative the the leading part. 

We use the notation $\theta = (\hat{\lambda},\omega)$ and apply the stationary phase lemma with respect to $\hat{t},\theta$ to compute the leading part as $|(\xi_{\mathrm{1c}},\eta_{\mathrm{1c}}) | \rightarrow \infty$.  We decompose $\theta$ according to directions parallel to and orthogonal to $(\xi_{\mathrm{1c}},\eta_{\mathrm{1c}})$ and denote projections of $\theta$ by $\theta^\parallel,\theta^\perp$ respectively. Then the critical set is given by $\hat{t}=0,\theta^\parallel = 0$. So the leading part, up to a constant factor, is
\begin{align}
|(\xi_{\mathrm{1c}},\eta_{\mathrm{1c}})|^{-1} x \int_{\mathbb{S}^{n-2}} \tilde{\chi}(x,y,\hat{\lambda}(\theta^\perp),\omega(\theta^\perp)) E_c d\theta^\perp, \label{ah_leading}
\end{align}
where the $|(\xi_{\mathrm{1c}},\eta_{\mathrm{1c}})|^{-1}$ comes from the square root of the determinant of the Hessian of the phase in the stationary phase lemma and $\hat{\lambda}(\theta^\perp),\omega(\theta^\perp)$ indicates that this critical set is parametrized by $\theta^\perp$ and thus other variables are functions of it. Recall (\ref{Ec_1form}), the tensorial factor on this critical set is
\begin{align*}
E_c = (\hat{\lambda}\frac{dx}{hx^3} + \frac{g_1(\omega\partial_y)}{h^{1/2}x}) \otimes( \hat{\lambda}hx^3\partial_x + \omega h^{1/2}x\partial_y).
\end{align*}
Since the statement is pointwise with respect to $x,y$, we assume that $g_1$ is the Euclidean metric in the ellipticity argument and $E_c$ becomes
\begin{align*}
E_c = (\hat{\lambda}\frac{dx}{hx^3} + \frac{\omega \cdot dy}{h^{1/2}x}) \otimes ( \hat{\lambda}hx^3\partial_x + \omega h^{1/2}x\partial_y).
\end{align*}

Since $\chi \geq 0$, this is a positive multiple of the projection to the span of $(\hat{\lambda},\omega)$. As $(\hat{\lambda},\omega)$ runs over the equatorial sphere consists of vectors orthogonal to $(\xi_{\mathrm{1c}},\eta_{\mathrm{1c}})$, we are integrating these projections, with the weight being strictly positive if $\chi(x,y,\hat{\lambda},\omega)>0$. 

Recall that the kernel of the principal symbol of $\delta_{h,\F}^s$ consists of covectors of the form $v=(v_0,v_1)$ such that $\xi_{\mathrm{1c}}v_0+\eta_{\mathrm{1c}}v'=0$. Thus we need to show that for each such $(v_0,v_1)$, there is at least one $(\hat{\lambda},\omega)$ in the critical set of the phase making $\tilde{\chi}>0$. Concretely, they satisfy

\begin{align}
& \tilde{\chi}(x,y,\hat{\lambda},\omega) > 0, \label{vcondition1}\\
& \xi_{\mathrm{1c}}\hat{\lambda} + \eta_{\mathrm{1c}} \cdot \omega = 0,\label{vcondition2}\\
& \hat{\lambda}v_0 + \omega \cdot v_1 \neq 0, \label{vcondition3}\\
& \xi_{\mathrm{1c}}v_0+\eta_{\mathrm{1c}}\cdot v_1 = 0. \label{vcondition4}
\end{align}

These conditions let us concludes that the integral of projections is positive and our desired ellipticity on the kernel of the standard principal symbol of $\delta_{h,\F}^s$ follows. So we arrange those conditions on vector components now. 

\begin{itemize}
\item Consider the case $v_1=0$ first. From (\ref{vcondition1}) we know $\xi_{\mathrm{1c}} = 0$. In order to satisfy (\ref{vcondition3}), we take $\hat{\lambda}$ small but non-zeri, $\omega$ orthogonal to $\eta_{\mathrm{1c}}$. Such $\omega$ exists because the orthogonal relationship is in $\R^{n-1},n \geq 3$.

\item Then we consider the case $v_1 \neq 0$ and $v_1$ is not parallel to $\eta$. In this case we take $\omega$ orthogonal to $\eta$ but not to $v_1$, $\hat{\lambda}=0$, then (\ref{vcondition1})-(\ref{vcondition3}) are satisfied.

\item Then we consider when $v_1 = c \eta \neq 0$, then (\ref{vcondition4}) becomes $\xi_{\mathrm{1c}}v_0+c|\eta_{\mathrm{1c}}|^2=0$. Then with $\omega$ to be determined later, we set $\hat{\lambda}=-\xi_{\mathrm{1c}}^{-1}\eta_{\mathrm{1c}}\cdot \omega$, we have
\begin{align*}
\hat{\lambda}v_0 + \omega \cdot v_1 = c(\omega \cdot \eta_{\mathrm{1c}})(1+ \xi_{\mathrm{1c}}^{-2}|\eta_{\mathrm{1c}}|^2),
\end{align*}
which is non-zero when $\omega \cdot \eta_{\mathrm{1c}} \neq 0$, which can be arranged since $\omega$ and $\eta_{\mathrm{1c}}$ are $(n-1)-$dimensional vectors. In addition, we may require $\omega \cdot \eta_{\mathrm{1c}}$ to be small so that $\hat{\lambda}$ is small as well and satisfy (\ref{vcondition1}).
\end{itemize}

The discussion on the contribution of the region $\{ |\hat{t}| \geq 1 \}$ is divided again into three parts: The first one is $\{|t| \leq T_0, \hat{t} \geq 1\}$ with $T_0$ fixed, the second one is the region $t$ bounded but away from 0 and the third region is $|t| \rightarrow \infty$. In the first region the phase is lower bounded by $|(\xi_{\mathrm{1c}},\eta_{\mathrm{1c}})| |\hat{t}|^{-k}$ for some $k$ and integration by parts gives the desired Schwartz property. In the second region, the no conjugate point assumption implies non-degenerate property of the Jacobian and in turn implies there is no critical point of the phase. The third region is analyzed in the similar manner to the second. See discussion after (3.13) of \cite{zachos2022inverting}, which verbatim transplant to our setting.

Combining cases above proves that the principal symbol of $N_{h,\F}$ is strictly positive definite on the kernel of principal symbol of $\delta_{h,\F}^s$.
\end{proof}

Next we consider the ellipticity at finite points.
\begin{lmn}
Let $\F>0$ and $\Omega_{x_0} = \{ x \leq x_0 \}$ with $x_0$ small.
$N_{h,\F}$ is an operator in
${\Psi}_{\mathrm{1c},h,\mathcal{F}}^{-1,-1}$ acting on one forms and it
is elliptic on the kernel of $\delta_{h,\F}$ in the boundary and semiclassical senses.
 \label{lmn_oneform_semiclassical}
\end{lmn}
\begin{proof}
Recall (\ref{symbol_oneform}), setting $h=0$ in the computation of semiclassical principal symbol, we have
\begin{align*}
\begin{split}
a_h(x,y,\lambda,\omega)= & x \int e^{\F( \hat{\lambda}\hat{t} +\alpha\hat{t}^2)} \tilde{\chi}(x,y,\lambda/h^{1/2},\omega)E_ce^{i(\xi_{\mathrm{1c}}(\hat{\lambda}\hat{t}+\alpha \hat{t}^2) + \eta_{\mathrm{1c}} \cdot  \omega \hat{t} )} d\hat{t}d\hat{\lambda} d\omega ,
\end{split}
\end{align*}
which is equivalent to
\begin{align}
\begin{split}
a_h(x,y,\lambda,\omega)= & x \int  e^{-i((-\xi_{\mathrm{1c}}+\F i)\hat{t})\hat{\lambda} }\tilde{\chi}(x,y,\lambda/h^{1/2},\omega)
\\& E_ce^{i\alpha \xi_{\mathrm{1c}} \hat{t}^2+i\eta_{\mathrm{1c}}\cdot\omega\hat{t}+\F\alpha\hat{t}^2}  d\hat{t}d\hat{\lambda} d\omega.
\end{split}
\end{align}
The $\hat{\lambda}-$integral is a Fourier transform evaluated at $(-\xi_{\mathrm{1c}}+\F i)$. Then multiplication by $\hat{\lambda}$ is transformed into $-D_\sigma=-i\partial_\sigma$, where $\sigma$ is the third variable of $\mathcal{F}_3\tilde{\chi}$ and $\mathcal{F}_3$ is the partial Fourier transform with respect to the third variable. Define the matrix $D_{c}$ to be $E_c$ with $\hat{\lambda}$ replaced by $-D_\sigma$. 
\begin{align*}
D_c =   \begin{pmatrix}
D_\sigma^2-\alpha\hat{t}D_\sigma && -D_\sigma \la \omega,\cdot \ra\\
\omega (-D_\sigma+2\alpha\hat{t}) && \omega \la \omega,\cdot \ra
\end{pmatrix}.
\end{align*}
This integral becomes ($\hat{\lambda}$ becomes $-D_\sigma$ after Fourier transform, because it becomes $D_\sigma$ after inverse Fourier transform)
\begin{align}
\begin{split}
a_h(x,y,\lambda,\omega)= & x \int e^{i\alpha \xi_{\mathrm{1c}} \hat{t}^2+i\eta_{\mathrm{1c}}\cdot\omega\hat{t}+\F\alpha\hat{t}^2} D_{c} \mathcal{F}_{3}\tilde{\chi}(x,y,((-\xi+i\F)\hat{t}),\omega) d\hat{t} d\omega ,
\end{split}
\end{align}

Take $\nu=\F^{-1}\alpha,\chi(s)=e^{\frac{s^2}{2 \nu } }$, then $\hat{\chi}(\sigma)=ce^{\frac{ \nu \sigma^2}{2}}$. Substitute in the equation above, we obtain  
\begin{align}
\begin{split}
a_h(x,y,\lambda,\omega)= &c x\int e^{i\alpha \xi_{\mathrm{1c}} \hat{t}^2+i\eta_{\mathrm{1c}} \cdot\omega\hat{t}+\F\alpha\hat{t}^2} D_{c} e^{\frac{\nu(-\xi_{\mathrm{1c}} +i\F)^2\hat{t}^2}{2}} d\hat{t} d\omega.
\end{split}
\end{align}

Notice that $D_\sigma$ is ($-i$ times) differentiating $(-\xi_{\mathrm{1c}}+i\F)\hat{t}$. For the convenience of later discussion, we set
\begin{align*}
\phi(\xi_{\mathrm{1c}},\omega) = & -\nu (-\xi_{\mathrm{1c}}+i\F)^2 - 2i\alpha\xi_{\mathrm{1c}}-2\F\alpha\\
                 = & -\nu (\xi_{\mathrm{1c}}^2+\F^2).
\end{align*} 
Then we compute derivatives of $\hat{\chi}=e^{\frac{\nu \sigma^2}{2}}$, with $c$ representing possibly different overall factors.
\begin{align*}
& D_\sigma \hat{\chi} = c \nu \sigma e^{\frac{\nu\sigma^2}{2}}, \quad  D^2_\sigma \hat{\chi} = c(\nu+\nu^2s^2) e^{\frac{\nu\sigma^2}{2}},
\end{align*}
Substituting those $\sigma-$derivatives back, we have
\begin{align*}
a_h(x,y,\hat{\lambda},\omega)= &c x \int_{\mathbb{S}^{n-2}} \int_{\R} e^{i \hat{t}\omega\cdot\eta_{\mathrm{1c}}}e^{-\frac{\phi\hat{t}^2}{2}}
\\&
\scriptsize
\begin{pmatrix}
-i\nu(\xi_{\mathrm{1c}}-i\F)(-i\nu(\xi_{\mathrm{1c}}-i\F)+2\alpha)\hat{t}^2-\nu
&& -i\nu(\xi_{\mathrm{1c}}-i\F)\hat{t}\la \omega,\cdot\ra \\   
\omega(-i\nu(\xi_{\mathrm{1c}}-i\F)+2\alpha)\hat{t} 
&& \omega \la \omega,\cdot \ra )
\end{pmatrix}  d\hat{t}d\omega.
\end{align*} 
The $\hat{t}$ integral is an inverse Fourier transform with respect to
$\omega \cdot \eta_{\mathrm{1c}}$, which turns multiplication by
$\hat{t}$ to $-D_{\omega \cdot \eta_{\mathrm{1c}}}$. In addition, $e^{-\frac{\phi\hat{t}^2}{2}}$ is transformed into (a constant multiple of) $\phi(\xi_{\mathrm{1c}},\omega)^{-1/2}e^{-\frac{(\omega\cdot \eta_{\mathrm{1c}})^2}{2\phi(\xi_{\mathrm{1c}},\omega)}}$. So we have
\begin{align*}
a_h(x,y,\hat{\lambda},\omega) &= cx \int_{\mathbb{S}^{n-2}} \int_{\R} \phi(\xi_{\mathrm{1c}},\omega)^{-1/2}e^{-\frac{(\omega\cdot \eta)^2}{2\phi}}\\
&
\scriptsize
\begin{pmatrix}
-i\nu(\xi_{\mathrm{1c}}-i\F)(-i\nu(\xi_{\mathrm{1c}}-i\F)+2\alpha)D_{\omega \cdot \eta_{\mathrm{1c}}}^2-\nu
&& -i\nu(\xi_{\mathrm{1c}}-i\F)D_{\omega \cdot \eta}\la \omega,\cdot\ra \\   
\omega(-i\nu(\xi_{\mathrm{1c}}-i\F)+2\alpha)D_{\omega \cdot \eta_{\mathrm{1c}}}
&& \omega \la \omega,\cdot \ra )
\end{pmatrix}d\omega.
\end{align*}
Direct computation shows this is 
\begin{align*} 
cx
(\xi_{\mathrm{1c}}^2+\F^2)^{-1/2}\int_{\mathbb{S}^{n-2}} & (-\nu)^{-1/2} 
\begin{pmatrix} 
 \frac{(\xi_{\mathrm{1c}}-i\F)}{\xi_{\mathrm{1c}}^2+\F^2}(\omega\cdot \eta_{\mathrm{1c}}) \\ 
\omega
\end{pmatrix} 
\\& \otimes \begin{pmatrix}  -\frac{(\xi_{\mathrm{1c}}+i\F)}{\xi_{\mathrm{1c}}^2+\F^2}(\omega \cdot \eta_{\mathrm{1c}}) && \la \omega,\cdot \ra      \end{pmatrix}
 e^{\frac{(\omega\cdot \eta_{\mathrm{1c}})^2}{2\nu(\xi_{\mathrm{1c}}^2+\F^2)}} d\omega.
\end{align*} 
We refer readers to page 21-22 of \cite{stefanov2018inverting} for a detailed process of a similar computation. We apply this principal symbol to an one form $v=v_0 \frac{dx}{h}+v_1\frac{dy}{h^{1/2}}$ in the kernel of the principal symbol of $\delta_h^s$, i.e. satisfying
\begin{align*}
(\xi_{\mathrm{1c}}+i\F)v_0 + \eta_{\mathrm{1c}} \cdot v_1 = 0.
\end{align*}
In order for the action of $a_{h}$ on $v$ to be non-vanishing, we need to choose $\omega$ so that
\begin{align}
-\frac{(\xi_{\mathrm{1c}}+i\F)}{\xi_{\mathrm{1c}}^2+\F^2}(\omega \cdot \eta_{\mathrm{1c}})v_0 + \omega \cdot v_1 = \frac{(\eta_{\mathrm{1c}}\cdot v_1)}{\xi_{\mathrm{1c}}^2+\F^2}(\omega \cdot \eta_{\mathrm{1c}}) + \omega \cdot v_1 \neq 0 .       \label{elliptic_1}
\end{align}
\begin{itemize}
\item If $v_1=0$, then $v=0$, which is trivial. \\
Next we consider the case in which $v_1 \neq 0$. 
\item If $\eta_{\mathrm{1c}}=0$, choosing $\omega$ parallel to $v_1$ gives (\ref{elliptic_1}). 
\item If $\eta_{\mathrm{1c}} \neq 0$, and $\eta_{\mathrm{1c}}$ is not parallel to $v_1$, then we take $\omega$ orthogonal to $\eta_{\mathrm{1c}}$ but not to $v_1$, which again gives (\ref{elliptic_1}). This is possible because $n \geq 3$ and $\omega$ has $n-2$ dimensional choices.
\item If $\eta \neq 0$ and $v_1 = c \eta$. Then the quantity in (\ref{elliptic_1}) is
\begin{align*}
c ( \frac{|\eta_{\mathrm{1c}}|^2}{\xi_{\mathrm{1c}}^2+\F^2} +1) (\omega \cdot v_1).       
\end{align*}
We can choose $\omega$ that is not orthogonal to $v_1$ to make (\ref{elliptic_1}) hold.
\end{itemize}
Summarizing all cases, $a_h$ is elliptic at finite points, and thus in
the boundary and semiclassical senses.
\end{proof}
\subsection{Ellipticity on 2-tensors}
We now consider the modified normal operator on 2-tensors.
Again, we first prove the ellipticity of $N_{h,\F}$ at fiber infinity.  
\begin{lmn}
Suppose $\F>0$ is sufficiently large and $\Omega_{x_0} = \{ x \leq x_0 \}$ with $x_0$ small. 
$N_{h,\F}$ is an operator in
${\Psi}_{\mathrm{1c},h,\mathcal{F}}^{-1,-1}$ acting on 2-tensors and it
is elliptic on the kernel of $\delta_{h,\F}$ in the 1-cusp algebra
differential sense.
 \label{lmn_2tensor_infinity}
\end{lmn}
\begin{proof}
Since the claim is pointwise, so we assume the $y-$part of $g_{h,\mathrm{1c}}$ is the $n-1$ dimensional Euclidean metric, i.e., $g_1=dy^2$. For 2-tensor computation, the basis is
\begin{align*}
\frac{dx}{hx^3} \otimes \frac{dx}{hx^3}, \, \frac{dx}{hx^3} \otimes \frac{dy}{h^{\frac{1}{2}}x}, \, \frac{dy}{h^{\frac{1}{2}}x} \otimes \frac{dx}{hx^3}, \, \frac{dy}{h^{\frac{1}{2}}x} \otimes \frac{dy}{h^{\frac{1}{2}}x} .
\end{align*} 

Similar to the $\tilde{E}_c$ factor in the one form case, we have an $\tilde{E}_{2c}$ factor in our 2-tensor case. In the 2-tensor case, the first factor becomes $g_{h,\mathrm{1c}}(\lambda x \partial_x + \omega \partial_y) \otimes g_{h,\mathrm{1c}}(\lambda x \partial_x + \omega \partial_y)$. For this factor, we have
\begin{align*}
 &g_{h,\mathrm{1c}}(\lambda x \partial_x + \omega \partial_y) \otimes g_{h,\mathrm{1c}}(\lambda x \partial_x + \omega \partial_y) 
\\= & h^{-1/2}x^{-1}(\hat{\lambda} \frac{dx}{hx^3}+\frac{g_1(\omega_1 \partial_y) }{h^{1/2}x}) \otimes h^{-1/2}x^{-1}(\hat{\lambda} \frac{dx}{hx^3}+\frac{g_1(\omega_2 \partial_y) }{h^{1/2}x})\\
= & h^{-1}x^{-2}\begin{pmatrix} \hat{\lambda}^2 \\ \hat{\lambda} \la  \omega,  \cdot \ra_1 \\ \hat{\lambda} \la  \omega,  \cdot \ra_2 \\ \la  \omega,  \cdot \ra_1 \la  \omega,  \cdot \ra_2  \end{pmatrix},
\end{align*}
where the indices $1,2$ in $\omega_1,\omega_2$ are indicating the order of components. While the second factor in the 2-tensor case is  
\begin{align*} 
\iota_{\dot{\gamma}} \otimes \iota_{\dot{\gamma}} = & h^{-1}x^{-2}( (\hat{\lambda}+2\alpha \hat{t})hx^3\partial_x + \omega_1 \partial_y ) \otimes ( (\hat{\lambda}+2\alpha \hat{t})hx^3\partial_x + \omega_2 \partial_y )\\
= & h^{-1}x^{-2} 
\setlength\arraycolsep{2pt}
\begin{pmatrix} (\hat{\lambda}+2\alpha\hat{t})^2 && (\hat{\lambda}+2\alpha \hat{t}) \la  \omega,  \cdot \ra_1 && (\hat{\lambda}+2\alpha \hat{t})\la  \omega,  \cdot \ra_2 &&  \la  \omega,  \cdot \ra_1 \la  \omega,  \cdot \ra_2   \end{pmatrix}.
\end{align*} 
$E_{2c}$, the product of two factors after combining $h$ and $x$ powers in the definition of $L$ and introduced by the change of variables is
\begin{align*}
\scriptsize
\begin{pmatrix}
\hat{\lambda}^2(\hat{\lambda}+2\alpha \hat{t})^2 && \hat{\lambda}^2(\hat{\lambda}+2\alpha\hat{t}) \la \omega,\cdot \ra_1 &&
\hat{\lambda}^2(\hat{\lambda}+2\alpha\hat{t})\la \omega,\cdot \ra_2 && \hat{\lambda}^2 \la \omega,\cdot \ra_1 \la \omega,\cdot \ra_2\\
 \hat{\lambda}(\hat{\lambda}+2\alpha\hat{t}) \omega_1 && \hat{\lambda}(\hat{\lambda}+2\alpha\hat{t}) \omega_1 \la \omega,\cdot \ra_1 && \hat{\lambda}(\hat{\lambda}+2\alpha \hat{t})\omega_1\la \omega,\cdot \ra_2 && \hat{\lambda}\omega_1 \la \omega,\cdot \ra_1 \la \omega,\cdot \ra_2\\
 \hat{\lambda}(\hat{\lambda}+2 \alpha \hat{t})^2\omega_2 && \hat{\lambda}(\hat{\lambda}+2\alpha \hat{t}) \omega_2 \la \omega, \cdot \ra_1 && \hat{\lambda}(\hat{\lambda}+2\alpha\hat{t})\omega_2 \la \omega,\cdot \ra_2 && \hat{\lambda} \omega_2 \la \omega,\cdot \ra_1 \la \omega,\cdot \ra_2 \\
 (\hat{\lambda}+2\alpha\hat{t})^2\omega_1\omega_2 && (\hat{\lambda}+2\alpha\hat{t})\omega_1\omega_2\la \omega,\cdot \ra_1 && (\hat{\lambda}+2\alpha\hat{t})\omega_1 \omega_2 \la \omega,\cdot \ra_2 && \omega_1\omega_2 \la \omega,\cdot \ra_1 \la \omega \cdot \ra_2
\end{pmatrix}.
\end{align*} 
The semiclassical principal symbol is
\begin{align}
\begin{split}
a_h(x,y,\lambda,\omega)= & x\int e^{\F( \hat{\lambda}\hat{t} +\alpha\hat{t}^2+\hat{t}^3h^{1/2}{\Gamma}^{(1)}(x,y,h^{1/2}\hat{\lambda},\omega,h^{1/2}\hat{t}) )} \tilde{\chi}(x,y,\lambda/xh^{1/2},\omega) E_{2c} \\
& e^{i(\xi(\hat{\lambda}\hat{t}+\alpha \hat{t}^2+h^{1/2}\hat{t}^2\Gamma^{(1)}) + \eta \cdot ( \omega \hat{t} + h^{1/2}\hat{t}^2\Gamma^{(2)}))} d\hat{t}d\hat{\lambda} d\omega ,
\end{split}
\end{align}
where $\hat{\lambda},\hat{t}$ are rescaled variables introduced in the one form case.
Following the computation in \cite{zachos2022inverting} and references therein, i.e., equation (3.12) of \cite{vasy2020semiclassical}, the symbol when we take $h^{1/2}=0$ is 
\begin{align*}
a_{h,2c}(x,y,\xi_{\mathrm{1c}},\eta_{\mathrm{1c}}) = x \int & e^{i\xi_{\mathrm{1c}}(\hat{\lambda}\hat{t}+\alpha(x,y,0,\omega)\hat{t}^2) + \eta_{\mathrm{1c}}\cdot \omega \hat{t} } e^{\F(\hat{\lambda}\hat{t} +\alpha(x,y,0,\omega)\hat{t}^2)} 
\\& \tilde{\chi}(x,y,\hat{\lambda},\omega)E_{2c}  d\hat{t}d\hat{\lambda}d\omega.
\end{align*}
Error terms caused by taking $h^{1/2}=0$ have extra $O(xh^{1/2} \la | (\xi_{\mathrm{1c}},\eta_{\mathrm{1c}}) | \ra^{-1} )$ gain relative the the leading part. 

We use the notation $\theta = (\hat{\lambda},\omega)$ and apply the stationary phase lemma with respect to $\hat{t},\theta$ to compute the leading part as $|(\xi_{\mathrm{1c}},\eta_{\mathrm{1c}}) | \rightarrow \infty$.  We decompose $\theta$ according to directions parallel to and orthogonal to $(\xi_{\mathrm{1c}},\eta_{\mathrm{1c}})$ and denote projections of $\theta$ by $\theta^\parallel,\theta^\perp$ respectively. Then the critical set is given by $\hat{t}=0,\theta^\parallel = 0$. So the leading part is
\begin{align*}
x \int_{\mathbb{S}^{n-2}} \tilde{\chi}(x,y,\hat{\lambda}(\theta^\perp),\omega(\theta^\perp))  E_{2c} d\theta^\perp,
\end{align*}
where $\hat{\lambda}(\theta^\perp),\omega(\theta^\perp)$ indicates that this critical set is parametrized by $\theta^\perp$ and thus other variables are functions of it.

Since $\chi \geq 0$, this is a positive multiple of the projection to the span of $(\hat{\lambda},\omega)$. As $(\hat{\lambda},\omega)$ runs over the equatorial sphere consists of vectors orthogonal to $(\xi_{\mathrm{1c}},\eta_{\mathrm{1c}})$, we are integrating these projections, with the weight being strictly positive if $\chi(x,y,\hat{\lambda},\omega)>0$. 

Recall that the kernel of the standard principal symbol of $\delta_{h,\F}^s$ consists of 2-tensors $v=(v_{00},v_{01},v_{01},v_{11})$ such that 
\begin{align}
\begin{split}
&\xi_{\mathrm{1c}} v_{00} + \eta_{\mathrm{1c}} \cdot v_{01}=0,\\
&\xi_{\mathrm{1c}} v_{01} + \frac{1}{2}(\eta_{\mathrm{1c},1}+\eta_{\mathrm{1c},2}) \cdot v_{11}=0,
\end{split}
\label{kernel_2tensor_infinity}
\end{align} 
where sub-indices in $\eta_{\mathrm{1c},1}$ resp. $\eta_{_{\mathrm{1c}},2}$ denotes the inner product is taken in the first resp. second slots of $v_{11}$. Taking the inner product with $\eta_{\mathrm{1c}}$ in the second equation, notice that $v_{11}$ is a 2-tensor sending a (co)vector $\eta_{\mathrm{1c}}$ to be a (co)vector, we obtain
\begin{align*}
\xi_{\mathrm{1c}} \eta_{\mathrm{1c}} \cdot \cdot v_{10} + (\eta_{\mathrm{1c}} \otimes \eta_{\mathrm{1c}}) \cdot v_{11} = 0.
\end{align*} 
Combining this with the first equation in (\ref{kernel_2tensor_infinity}) yields
\begin{align*}
\xi_{\mathrm{1c}}^2 v_{00}  = (\eta_{\mathrm{1c}} \otimes \eta_{\mathrm{1c}})v_{11}.
\end{align*} 

Without $i\F-$term in the finite points case, we need to consider cases in which $\xi=0$ and $\xi \neq 0$ respectively.

Consider the case $\xi \neq 0$ first.
\begin{align}
\begin{split}
&v_{00}=\xi_{\mathrm{1c}}^{-2}(\eta\otimes\eta)v_{11},\\
&v_{01} = -\frac{1}{2\xi_{\mathrm{1c}}} (\eta_{\mathrm{1c},1}+\eta_{\mathrm{1c},2})\cdot v_{11}.
\end{split}
\label{2c_v}
\end{align} 

Recall that we are computing the contribution at the critical set of the phase: $\{\hat{t}=0,\theta^\parallel = 0\}$, so the direction to which we are projecting is

\begin{align*} 
\iota_{\dot{\gamma}} \otimes \iota_{\dot{\gamma}} = & h^{-1}x^{-2} \begin{pmatrix} ( \hat{\lambda}^2 && \hat{\lambda} \la  \omega,  \cdot \ra_1 && \hat{\lambda}\la  \omega,  \cdot \ra_2 &&  \la  \omega,  \cdot \ra_1 \la  \omega,  \cdot \ra_2   \end{pmatrix}.
\end{align*} 

As $\theta^\parallel=0$, or equivalently $\xi_{\mathrm{1c}} \hat{\lambda}  + \eta_{\mathrm{1c}} \omega = 0$, we know
\begin{align*}
\hat{\lambda} = - \frac{\eta_{\mathrm{1c}} \cdot \omega}{\xi_{\mathrm{1c}}}.
\end{align*}

So for $v$ to be in the kernel of the projection means
\begin{align*}
( (\frac{\eta_{\mathrm{1c}} \cdot \omega}{\xi_{\mathrm{1c}}})^2 \frac{\eta_{\mathrm{1c}} \otimes \eta_{\mathrm{1c}}}{\xi_{\mathrm{1c}}^2}   + 
\frac{\eta_{\mathrm{1c}} \cdot \omega}{\xi_{\mathrm{1c}}^2}( \frac{\eta_{\mathrm{1c}}}{\xi_{\mathrm{1c}}} \otimes \omega + \omega \otimes \frac{\eta_{\mathrm{1c}}}{\xi_{\mathrm{1c}}})  + \omega \otimes \omega) \cdot v_{11} = 0,
\end{align*} 

which is equivalent to
\begin{align}
( ( \frac{\eta_{\mathrm{1c}} \cdot \omega}{\xi_{\mathrm{1c}}} \frac{\eta_{\mathrm{1c}}}{\xi_{\mathrm{1c}}} +\omega) \otimes ( \frac{\eta_{\mathrm{1c}} \cdot \omega}{\xi_{\mathrm{1c}}} \frac{\eta_{\mathrm{1c}}}{\xi_{\mathrm{1c}}} +\omega) ) \cdot v_{11} = 0. 
\label{2c_v2}
\end{align}

Suppose $\eta_{\mathrm{1c}} = 0$, then (\ref{2c_v}) implies $v_{00}=0,v_{01}=0$. Equation (\ref{2c_v2}) now becomes
\begin{align*}	
(\omega \otimes \omega) \cdot v_{11} = 0.
\end{align*}
By the `polarization' formula $\omega_1 \otimes \omega_2+\omega_2\otimes \omega_1=(\omega_1+\omega_2)\otimes(\omega_1+\omega_2)-\omega_1\otimes\omega_1-\omega_2\otimes\omega_2$, we know that the symmetric 2-tensors of the form $\omega \otimes \omega$ span the sapce of all symmetric 2-tensors. So from $(\omega \otimes \omega) \cdot v_{11} = 0$ we conclude that $v_{11}=0$ and consequently $v=0$ in this case.

Now suppose $\eta_{\mathrm{1c}} \neq 0$, with $\hat{\eta}_{\mathrm{1c}}=\frac{\eta_{\mathrm{1c}}}{|\eta_{\mathrm{1c}}|}$, we decompose $\omega$ as 
\begin{align*}
\omega = \epsilon \hat{\eta}_{\mathrm{1c}}+(1-\epsilon^2)^{1/2}\omega^\perp,
\end{align*}
 where $\omega^\perp$ is the unit vector in the direction of the projection of $\omega$ onto $\eta_{\mathrm{1c}}^\perp$, the orthogonal complement of $\eta_{\mathrm{1c}}$. Then we substitute this into (\ref{2c_v2}) to obtain
\begin{align*}
& ( (1+ \frac{|\eta_{\mathrm{1c}}|^2}{\xi_{\mathrm{1c}}^2})^2\epsilon^2 \hat{\eta_{\mathrm{1c}}} \otimes \hat{\eta}_{\mathrm{1c}})+ (1+ \frac{|\eta_{\mathrm{1c}}|^2}{\xi_{\mathrm{1c}}^2})\epsilon(1-\epsilon^2)^{1/2}(\hat{\eta}_{\mathrm{1c}} \otimes \omega^{\perp} +\omega^\perp \otimes \hat{\eta}_{\mathrm{1c}} )\\
& +(1-\epsilon^2)\omega^\perp\otimes\omega^\perp)\cdot v_{11}=0.
\end{align*}
Let the directions of $\omega$ vary, then $\epsilon$ varies correspondingly, and our condition is that this equation holds for all $\epsilon \in [-1,1]$. In addition, now $\eta_{\mathrm{1c}} \cdot \omega = \epsilon |\eta_{\mathrm{1c}}|$ and we have $\hat{\lambda} = -\frac{\epsilon |\eta_{\mathrm{1c}}|}{\xi_{\mathrm{1c}}}$, so $\hat{\lambda}$ is small when $\epsilon$ is small and therefore $\tilde{\chi}>$ for small $\epsilon$. Taking $\epsilon = 0$ first yields $(\omega^\perp \otimes \omega^\perp) \cdot v_{11} = 0$.
By the same polarization argument as above, cotensors of the form $\omega^\perp \otimes \omega^\perp$ span $\eta_{\mathrm{1c}} \otimes \eta_{\mathrm{1c}}^\perp$, we conclude that $v_{11}$ is orthogonal to every cotensor in $\eta_{\mathrm{1c}}^\perp \otimes \eta_{\mathrm{1c}}^\perp$. 
Our second step is to take derivative with respect to $\epsilon$ at $\epsilon=0$ at (\ref{kernel_decomp}), which yields $(\hat{\eta}_{\mathrm{1c}} \otimes \omega^\perp + \omega^\perp \otimes \hat{\eta}_{\mathrm{1c}}) \cdot v_{11}=0$ for all $\omega^\perp$. Notice that symmetric tensors of the form $(\hat{\eta}_{\mathrm{1c}} \otimes \omega^\perp + \omega^\perp \otimes \hat{\eta}_{\mathrm{1c}})$ and $\eta_{\mathrm{1c}}^\perp \otimes \eta_{\mathrm{1c}}^\perp$ together span $(\eta_{\mathrm{1c}} \otimes \eta_{\mathrm{1c}})^\perp$. Then further taking the second order derivative at $\epsilon=0$ shows that $(\hat{\eta}_{\mathrm{1c}} \otimes \hat{\eta}_{\mathrm{1c}}) \cdot v_{11}=0$. Combining orthogonal conditions we have we know $v_{11}=0$. Combining (\ref{2c_v}), we know that $v=0$, then the non-degeneracy of this principal symbol and thus the ellipticity when $\xi_{\mathrm{1c}} \neq 0$ follows.

Now we consider the case $\xi_{\mathrm{1c}} = 0$. Now the critical set condition becomes $\eta_{\mathrm{1c}} \cdot \omega =0$. (\ref{2c_v}) and the condition being in the kernel of $\iota_{\dot{\gamma}}\otimes \iota_{\dot{\gamma}}$ is
\begin{align}
\begin{split}
& \eta_{\mathrm{1c}} \cdot v_{01} = 0,\\
& (\eta_{\mathrm{1c},1}+\eta_{\mathrm{1c},2}) \cdot v_{11} = 0,\\
& ( \hat{\lambda}^2 v_{00}+2\hat{\lambda}\omega \cdot v_{01} + (\omega \otimes \omega) \cdot v_{11}) = 0.
\end{split} \label{xi=0}
\end{align}
Suppose (\ref{ah_leading}) acting on $v$ is not elliptic, then it vanishes for every $\hat{\lambda},\omega$ satisfy $\eta_{\mathrm{1c}} \cdot \omega =0$ and $\hat{\lambda}$ in the support of $\tilde{\chi}$ (for fixed other variables). In particular, for every $\hat{\lambda}$ small. 
View this as a polynomial in $\hat{\lambda}$, then its coefficients need to vanish. So we conclude that $v_{00}=0, \omega \cdot v_{01} = 0 , (\omega \otimes \omega) \cdot v_{11} = 0 $. Combine these with first two equations in (\ref{xi=0}), we know $v_{01}$ is orthogonal to both $\eta_{\mathrm{1c}}$ and all directions orthogonal to it, similarly for $v_{11}$ (using the argument in the previous $\xi \neq 0$ part), we conclude that $v_{01} = 0, v_{11}=0$, thus $v=0$ as desired. The proof of ellipticity at fiber infinity is completed.
\end{proof}

Next we consider the ellipticity at finite points.
\begin{lmn}
Suppose $\F>0$ is sufficiently large and $\Omega_{x_0} = \{ x \leq x_0 \}$ with $x_0$ small.
$N_{h,\F}$ is an operator in
${\Psi}_{\mathrm{1c},h,\mathcal{F}}^{-1,-1}$ acting on 2-tensors and it
is elliptic on the kernel of $\delta_{h,\F}$ in the boundary and semiclassical senses.
 \label{lmn_2tensor_semiclassical}
\end{lmn}
\begin{proof}
Then the semiclassical principal symbol at finite points is
\begin{align*}
a_{h,2c}:=x \int e^{ \alpha \hat{t}^2 (\F+i\xi_{\mathrm{1c}}) + \hat{t}(\hat{\lambda}(\F+i\xi_{\mathrm{1c}})+i\eta_{\mathrm{1c}}\cdot \omega) } \tilde{\chi}(x,y,\hat{\lambda},\omega) 
E_{2c} d\hat{t} d\hat{\lambda} d\omega.
\end{align*} 
Use the expression from \cite{zachos2022inverting}, adding the tensorial factor, and notice that with our new parameter $\F$, the phase is unchanged, while the damping factor is multiplied by $\F$, we have
\begin{align}
\begin{split}
a_h(x,y,\lambda,\omega)= & x \int e^{\F( \hat{\lambda}\hat{t} +\alpha\hat{t}^2)} \tilde{\chi}(x,y,\lambda/h^{1/2},\omega)E_{2c}e^{i(\xi_{\mathrm{1c}}(\hat{\lambda}\hat{t}+\alpha \hat{t}^2) + \eta_{\mathrm{1c}} \cdot  \omega \hat{t} )} d\hat{t}d\hat{\lambda} d\omega ,
\end{split}
\end{align}
Rewrite $a_h$ as
\begin{align}
\begin{split}
a_h(x,y,\lambda,\omega)=  x\int &  e^{-i((-\xi_{\mathrm{1c}}+\F i)\hat{t})\hat{\lambda} }\tilde{\chi}(x,y,\lambda/h^{1/2},\omega)
\\ &{E}_{2c}e^{i\alpha \xi_{\mathrm{1c}} \hat{t}^2+i\eta_{\mathrm{1c}}\cdot\omega\hat{t}+\F\alpha\hat{t}^2}  d\hat{t}d\hat{\lambda} d\omega .
\end{split}
\label{symbol_1c_E2c}
\end{align}
The $\hat{\lambda}-$integral is a Fourier transform evaluated at $(-\xi_{\mathrm{1c}}+\F i)$. Then multiplication by $\hat{\lambda}$ is transformed into $-D_\sigma=i\partial_\sigma$, where $\sigma$ is the third variable of $\mathcal{F}_3\tilde{\chi}$. Define the matrix $D_{2c}$ to be $E_{2c}$ with $\hat{\lambda}$ replaced by $-D_\sigma$:
\begin{align*}
\scriptsize
\setlength\arraycolsep{2pt}
\begin{pmatrix}
D_\sigma^2(-D_\sigma+2\alpha \hat{t})^2 & D_\sigma^2(-D_\sigma+2\alpha\hat{t}) \la \omega,\cdot \ra_1 &
D_\sigma^2(-D_\sigma+2\alpha\hat{t})\la \omega,\cdot \ra_2 & D_\sigma^2 \la \omega,\cdot \ra_1 \la \omega,\cdot \ra_2\\
 -D_\sigma(-D_\sigma+2\alpha\hat{t}) \omega_1 & -D_\sigma(-D_\sigma+2\alpha\hat{t}) \omega_1 \la \omega,\cdot \ra_1 & -D_\sigma(-D_\sigma+2\alpha \hat{t})\omega_1\la \omega,\cdot \ra_2 & -D_\sigma\omega_1 \la \omega,\cdot \ra_1 \la \omega,\cdot \ra_2\\
 -D_\sigma(-D_\sigma+2 \alpha \hat{t})^2\omega_2 & -D_\sigma(-D_\sigma+2\alpha \hat{t}) \omega_2 \la \omega, \cdot \ra_1 & -D_\sigma(-D_\sigma+2\alpha\hat{t})\omega_2 \la \omega,\cdot \ra_2 & -D_\sigma \omega_2 \la \omega,\cdot \ra_1 \la \omega,\cdot \ra_2 \\
 (-D_\sigma+2\alpha\hat{t})^2\omega_1\omega_2 & (-D_\sigma+2\alpha\hat{t})\omega_1\omega_2\la \omega,\cdot \ra_1 & (-D_\sigma+2\alpha\hat{t})\omega_1 \omega_2 \la \omega,\cdot \ra_2 & \omega_1\omega_2 \la \omega,\cdot \ra_1 \la \omega \cdot \ra_2
\end{pmatrix}.
\end{align*} 
Then (\ref{symbol_1c_E2c}) becomes
\begin{align}
\begin{split}
a_h(x,y,\lambda,\omega)= & x \int e^{i\alpha \xi_{\mathrm{1c}} \hat{t}^2+i\eta_{\mathrm{1c}}\cdot\omega\hat{t}+\F\alpha\hat{t}^2} D_{2c} \mathcal{F}_{3}\tilde{\chi}(x,y,((-\xi_{\mathrm{1c}}+i\F)\hat{t}),\omega) d\hat{t} d\omega ,
\end{split}
\label{symbol_2c_f3}
\end{align}
Take $\nu=\F^{-1}\alpha,\tilde{\chi}(x,y,s,\omega)=e^{\frac{s^2}{2\nu}}$, where the dependence on $x,y,\omega$ is encoded in $\alpha$ and thus in $\nu$, 
then $\mathcal{F}_3\tilde{\chi}(x,y,\sigma,\omega)=ce^{\frac{ \nu \sigma^2}{2}}$. Substitute in (\ref{symbol_2c_f3}), we obtain  
\begin{align}
\begin{split}
a_h(x,y,\lambda,\omega)= &c x \int e^{i\alpha \xi_{\mathrm{1c}} \hat{t}^2+i\eta_{\mathrm{1c}} \cdot\omega\hat{t}+\F\alpha\hat{t}^2} D_{2c} e^{\frac{\nu(-\xi_{\mathrm{1c}}+i\F)^2\hat{t}^2}{2}} d\hat{t} d\omega.
\end{split}
\end{align} 
For the convenience of later discussion, we set
\begin{align*}
\phi(\xi_{\mathrm{1c}},\omega) = & -\nu (-\xi_{\mathrm{1c}}+i\F)^2 - 2i\alpha\xi_{\mathrm{1c}}-2\F\alpha\\
                 = & -\nu (\xi_{\mathrm{1c}}^2+\F^2).
\end{align*} 
Then we compute derivatives of $\mathcal{F}_3\tilde{\chi}=e^{\frac{\nu \sigma^2}{2}}$, with $c$ in in different equations representing possibly different overall factors.
\begin{align*}
& D_\sigma \hat{\chi} = c \nu \sigma e^{\frac{\nu\sigma^2}{2}}, \quad  D^2_\sigma \hat{\chi} = c(\nu+\nu^2s^2) e^{\frac{\nu\sigma^2}{2}},\\
& D_\sigma^3 \hat{\chi} = c(3\nu^2\sigma+\nu^3\sigma^3)e^{\frac{\nu\sigma^2}{2}} , \quad D_\sigma^4 \hat{\chi} = c(3\nu^2+6\nu^3\sigma^2+\nu^4\sigma^4)e^{\frac{\nu\sigma^2}{2}}.
\end{align*}
Then we have
\begin{align*}
a_h(x,y,\hat{\lambda},\omega)= &cxh \int_{\mathbb{S}^{n-2}} \int_{\R} e^{i \hat{t}\omega\cdot\eta_{\mathrm{1c}}}
(\tilde{B}_{ij})\times e^{-\frac{\phi\hat{t}^2}{2}}d\hat{t}d\omega.
\end{align*} 

Let $B_{ij}$ be the coefficient after the action of $D_\sigma^i(D_\sigma+2\alpha\hat{t})^j$, we have (notice that the variable of $\hat{\chi}$ is $-\xi_{\mathrm{1c}}+i\F$, but now we are writting expressions for $(\xi_{\mathrm{1c}}-i\F)$ for convenience)
\begin{align*}
& B_{00}=1,\\
& B_{10}=-i\nu(\xi_{\mathrm{1c}}-i\F)\hat{t},\\
& B_{20}=-\nu^2(\xi_{\mathrm{1c}}-i\F)^2\hat{t}^2-\nu,\\
& B_{01}=-i\nu(\xi_{\mathrm{1c}}-i\F)\hat{t}+2\alpha\hat{t}, \\
& B_{11}=- \nu(\xi_{\mathrm{1c}}-i\F)( \nu(\xi_{\mathrm{1c}}-i\F)+2i\alpha )\hat{t}^2 - \nu,\\
& B_{21}=i\nu^2(\xi_{\mathrm{1c}}-i\F)^2(\nu(\xi_{\mathrm{1c}}-i\F) + 2i\alpha)\hat{t}^3 + (3i\nu^2(\xi_{\mathrm{1c}}-i\F)-2\nu\alpha)\hat{t},\\
&B_{02}=(\nu(\xi_{\mathrm{1c}}-i\F)+2i\alpha)^2\hat{t}^2-\nu,\\
& B_{12}= -i\nu(\xi_{\mathrm{1c}}-i\F)(\nu(\xi_{\mathrm{1c}}-i\F) + 2i\alpha )^2\hat{t}^3 ,\\
& B_{22}=\nu^2(\xi_{\mathrm{1c}}-i\F)^2(\nu(\xi_{\mathrm{1c}}-i\F)+2i\alpha)^2\hat{t}^4+\nu(6\nu^2(\xi_{\mathrm{1c}}-i\F)^2
\\& \quad \quad +12i\nu\alpha(\xi_{\mathrm{1c}}-i\F)-4\nu^2)\hat{t}^2+3\nu^2.
\end{align*}
Then the matrix $(\tilde{B}_{ij})$ is given by
\begin{align*}
\begin{pmatrix}
B_{22} && B_{21} \la \omega,\cdot \ra_1 && B_{21} \la \omega,\cdot \ra_2 && B_{20}\la \omega,\cdot \ra_1 \la \omega,\cdot \ra_2\\
B_{12}\omega_1 && B_{11}\omega_1\la \omega,\cdot \ra_1 && B_{11}\omega_1\la \omega,\cdot \ra_2 && B_{10}\omega_1 \la \omega_1,\cdot \ra_1 \la \omega \ra_2 \\
B_{12}\omega_2 && B_{11}\omega_2 \la \omega,\cdot \ra_1 && B_{11}\omega_2 \la \omega,\cdot \ra_2 && B_{10}\omega_2 \la \omega_1,\cdot \ra_1 \la \omega \ra_2 \\
B_{02}\omega_1\omega_2 && B_{01}\omega_1\omega_2 \la \omega,\cdot \ra_1 &&  B_{01}\omega_1\omega_2 \la \omega,\cdot \ra_2 && B_{10}\omega_2\la \omega_2 \la \omega,\cdot \ra_1 \\
\end{pmatrix}
\end{align*}

The $\hat{t}$ integral is an inverse Fourier transform with respect to
$\omega \cdot \eta_{\mathrm{1c}}$, which turns multiplication by
$\hat{t}$ to $-D_{\omega \cdot \eta_{\mathrm{1c}}}$. Moreover, $e^{-\frac{\phi\hat{t}^2}{2}}$ is transformed into (a constant multiple of) $\phi(\xi_{\mathrm{1c}},\omega)^{-1/2}e^{-\frac{(\omega\cdot \eta_{\mathrm{1c}})^2}{2\phi(\xi_{\mathrm{1c}},\omega)}}$.
 So we have
\begin{align*}
a_h(x,y,\hat{\lambda},\omega)= &cx \int_{\mathbb{S}^{n-2}} \int_{\R} \phi(\xi_{\mathrm{1c}},\omega)^{-1/2}
(\tilde{C}_{ij}) \times e^{-\frac{(\omega\cdot \eta_{\mathrm{1c}})^2}{2\phi}}d\omega.
\end{align*}

The matrix $(\tilde{C}_{ij})$ is given by
\begin{align*}
\begin{pmatrix}
C_{22} && C_{21} \la \omega,\cdot \ra_1 && C_{21} \la \omega,\cdot \ra_2 && C_{20}\la \omega,\cdot \ra_1 \la \omega,\cdot \ra_2\\
C_{12}\omega_1 && C_{11}\omega_1\la \omega,\cdot \ra_1 && C_{11}\omega_1\la \omega,\cdot \ra_2 && C_{10}\omega_1 \la \omega_1,\cdot \ra_1 \la \omega \ra_2 \\
C_{12}\omega_2 && C_{11}\omega_2 \la \omega,\cdot \ra_1 && C_{11}\omega_2 \la \omega,\cdot \ra_2 && C_{10}\omega_2 \la \omega_1,\cdot \ra_1 \la \omega \ra_2 \\
C_{02}\omega_1\omega_2 && C_{01}\omega_1\omega_2 \la \omega,\cdot \ra_1 &&  C_{01}\omega_1\omega_2 \la \omega,\cdot \ra_2 && C_{10}\omega_2\la \omega_2 \la \omega,\cdot \ra_1 \\
\end{pmatrix}.
\end{align*}
Here $(\tilde{C}_{ij})$ is the counterpart of (3.20) of \cite{stefanov2018inverting}, set $\rho=\omega\cdot\eta_{\mathrm{1c}}$,
\begin{align*}
& C_{00}=1,\\
& C_{10}=\nu(\xi_{\mathrm{1c}}-i\F)\phi^{-1}\rho,\\
& C_{20}=\nu^2(\xi_{\mathrm{1c}}-i\F)^2\phi^{-2}\rho^2+2i\nu\alpha\phi^{-1}(\xi_{\mathrm{1c}}-i\F),\\
& C_{01}=\nu(\xi_{\mathrm{1c}}-i\F)\phi^{-1}\rho-2i\alpha\phi^{-1}\rho, \\
& C_{11}=\nu(\xi_{\mathrm{1c}}-i\F)( \nu(\xi_{\mathrm{1c}}-i\F)+2i\alpha )\phi^{-2}\rho^2,\\
& C_{21}=i\nu^2(\xi_{\mathrm{1c}}-i\F)^2(\nu(\xi_{\mathrm{1c}}-i\F) + 2i\alpha)\phi^{-3}\rho^3 - 2i\nu\alpha\phi^{-1}\rho,\\
& C_{02}=(\nu(\xi_{\mathrm{1c}}-i\F)+2i\nu)^2\phi^{-2}\rho^2-\phi^{-1}(\nu(\xi-i\F)+2i\nu)2i\nu,\\
& C_{12}= \nu(\xi_{\mathrm{1c}}-i\F)(\nu(\xi_{\mathrm{1c}}-i\F) + 2i\alpha )^2\phi^{-3}\rho^3+2i\nu^2\phi^{-1}\rho, \\
& C_{22}=\nu^2(\xi_{\mathrm{1c}}-i\F)^2(\nu(\xi_{\mathrm{1c}}-i\F)+2i\alpha)^2\phi^{-4}\rho^4+4\nu\alpha^2\rho^{-2}\phi^2-4\alpha^2\nu\phi^{-1}.\\
\end{align*}
We can verify that $C_{i0}C_{0j}=C_{ij}$, which is as expected since we constructed $E_{2c}$ and therefore subsequent matrices by multiplying two rank one matrices. Thus the matrix $(\tilde{C}_{ij})$ can be decomposed as
\begin{align}
\begin{pmatrix} C_{20} \\ \omega_1C_{10} \\ \omega_2 C_{10} \\ \omega_1\omega_2 \end{pmatrix} \otimes \begin{pmatrix} C_{02} && C_{01}\la \omega, \cdot \ra_1 && C_{01}\la \omega, \cdot \ra_2 && \la \omega, \cdot \ra_1 \la \omega, \cdot \ra_2 \end{pmatrix}  ,         \label{decomp_principal}
\end{align}
where the second factor is the adjoint of the first one. In addition, we have
\begin{align*}
& C_{01}=\nu(\xi_{\mathrm{1c}}+i\F)\phi^{-1}\rho,\\
& C_{02}=\nu^2(\xi_{\mathrm{1c}}+i\F)^2\phi^{-2}\rho^2+2i\nu\alpha\phi^{-1}(\xi_{\mathrm{1c}}+i\F),
\end{align*}
where $\phi=-\nu(\xi^2+\F^2)$. 

Using (\ref{deltas_symbol}), the condition that a symmetric 2-tensor $v=(v_{00},v_{01},v_{10},v_{11})$ (being symmetric means $v_{01}=v_{10}$ being in the kernel of the principal symbol of $\delta_h^s$ means that
\begin{align}
\begin{split}
&(\xi_{\mathrm{1c}}+i\F) v_{00} + \eta_{\mathrm{1c}} \cdot v_{01}+b_s\cdot v_{11}=0,\\
&(\xi_{\mathrm{1c}}+i\F) v_{01} + \frac{1}{2}(\eta_{\mathrm{1c},1}+\eta_{\mathrm{1c},2})\cdot v_{11}=0,
\end{split}
\label{kernel_2tensor_finite}
\end{align} 
where sub-indices in $\eta_{\mathrm{1c},1}$ resp. $\eta_{\mathrm{1c},2}$ denotes the inner product is taken in the first resp. second slots of $v_{11}$. Taking the inner product with $\eta_{\mathrm{1c}}$ in the second equation, notice that $v_{11}$ is a 2-tensor sending a (co)vector $\eta_{\mathrm{1c}}$ to be a (co)vector, we obtain
\begin{align*}
(\xi_{\mathrm{1c}}+i\F) \eta_{\mathrm{1c}} \cdot  v_{01} + (\eta_{\mathrm{1c}} \otimes \eta_{\mathrm{1c}}) \cdot v_{11} = 0.
\end{align*} 
Combining this with the first equation in (\ref{kernel_2tensor_finite}) yields
\begin{align*}
(\xi_{\mathrm{1c}}+i\F)^2 v_{00} + ((\xi_{\mathrm{1c}}+i\F)b_s  - \eta_{\mathrm{1c}} \otimes \eta_{\mathrm{1c}} ) \cdot v_{11}=0.
\end{align*}
From the second equation of (\ref{kernel_2tensor_finite}), we know
\begin{align*}
(\xi_{\mathrm{1c}}+i\F)v_{01} =-\frac{1}{2}(\eta_{\mathrm{1c},1}+\eta_{\mathrm{1c},2})\cdot v_{11}.
\end{align*}
Using symmetry of tensors involved to combine inner products taken, we have
\begin{align}
\begin{split}
&v_{00}=(\xi_{\mathrm{1c}}+i\F)^{-2}(\eta_{\mathrm{1c}}\otimes\eta_{\mathrm{1c}}- (\xi_{\mathrm{1c}}+i\F)b_s) \cdot v_{11},\\
&v_{01} = -\frac{1}{2}(\xi_{\mathrm{1c}}+i\F)^{-1}(\eta_{\mathrm{1c},1}+\eta_{\mathrm{1c},2})\cdot v_{11}.
\end{split}
\label{v_components_relation}
\end{align} 
Fixing $\omega$, for $v$ to be in the kernel of the principal symbol of $\delta_h^s$ and the projection given by (\ref{decomp_principal}) means
\begin{align*}
&(C_{02}(\xi_{\mathrm{1c}}+i\F)^{-2}(\eta_{\mathrm{1c}}\otimes\eta_{\mathrm{1c}}-(\xi_{\mathrm{1c}}+i\F)b_s) 
\\&- C_{01}(\xi_{\mathrm{1c}}+i\F)^{-1}(\eta_{\mathrm{1c}} \otimes \omega+\omega\otimes \eta_{\mathrm{1c}})+ \omega \otimes \omega ) \cdot v_{11} = 0.
\end{align*} 
Recalling that $\phi=-\nu(\xi_{\mathrm{1c}}^2+\F^2)$, we conclude from
the concrete expressions of $C_{ij}$ that this is equivalent to
\begin{align*}
& ((\xi_{\mathrm{1c}}-i\F)^{-1}(\xi_{\mathrm{1c}}^2+\F^2)^{-1}(\omega\cdot\eta_{\mathrm{1c}})^2 -2i\nu(\xi_{\mathrm{1c}}^2+\F^2)^{-1}) 
\\&((\xi_{\mathrm{1c}}+i\F)^{-1}(\eta_{\mathrm{1c}} \otimes \eta_{\mathrm{1c}}) -b_s) + 
\\& (\xi_{\mathrm{1c}}^2+\F^2)^{-1}(\omega\cdot\eta_{\mathrm{1c}})(\eta_{\mathrm{1c}} \otimes \omega+ \omega \otimes \eta_{\mathrm{1c}})+\omega \otimes \omega) \cdot v_{11} = 0.
\end{align*} 
In order to eliminate `error terms' involving $b_s$ above, we introduce 
\begin{align*}
\xi_{\mathrm{1c},\F} = \frac{\xi_{\mathrm{1c}}}{\F}  , \, \eta_{\mathrm{1c},\F} = \frac{\eta_{\mathrm{1c}}}{\F}.
\end{align*}
Using these rescaled variables, condition above can be rewritten as
\begin{align*}
&((\xi_{\mathrm{1c},\F}-i)^{-1}(\xi_{\mathrm{1c},\F}^2+1)^{-1}(\omega \cdot \eta_{\mathrm{1c},\F})^2 + 2i\F^{-1}\alpha(\xi_{\mathrm{1c},\F}^2+1)^{-1}) 
\\&((\xi_{\mathrm{1c},\F}+i)^{-1}(\eta_{\mathrm{1c},\F} \otimes \eta_{\mathrm{1c},\F})-\F^{-1}a) + \\
&(\xi_{\mathrm{1c},\F}^2+1)^{-1}(\omega \cdot \eta_{\mathrm{1c},\F})( \eta_{\mathrm{1c},\F} \otimes \omega + \omega \otimes \eta_{\mathrm{1c},\F} ) + \omega \otimes \omega ) \cdot v_{11}=0.
\end{align*}

We collect terms involving $\F^{-1}$ to rewrite this as:
\begin{align*}
& ( (\xi_{\mathrm{1c},\F}-i)^{-1}(\xi_{\mathrm{1c},\F}^2+1)^{-1}(\omega\cdot \eta_{\mathrm{1c},\F})^2)(\xi_{\mathrm{1c},\F}+i)^{-1}(\eta_{\mathrm{1c},\F} \otimes \eta_{\mathrm{1c},\F}) +\\
& (\xi_{\mathrm{1c},\F}^2+1)^{-1}(\omega \cdot \eta_{\mathrm{1c},\F})(\eta_{\mathrm{1c},\F} \otimes \omega + \omega \otimes \eta_{\mathrm{1c},\F}) + \omega \otimes \omega +O(\F^{-1})) \cdot v_{11} = 0.
\end{align*}
This can be decomposed as
\begin{align}
\begin{split}
&(( (\xi_{\mathrm{1c},\F}^2+1)^{-1}(\omega \cdot \eta_{\mathrm{1c},\F})\eta_{\mathrm{1c},\F} +\omega)
\\&\otimes( (\xi_{\mathrm{1c},\F}^2+1)^{-1}(\omega \cdot \eta_{\mathrm{1c},\F})\eta_{\mathrm{1c},\F} +\omega ) +O(\F^{-1}))\cdot v_{11}=0.
\end{split}
\label{kernel_decomp} 
\end{align}
Then we choose $\F$ large, and we show that this equation holds for all $\omega$ implies $v_{11}=0$. When $\eta_{\mathrm{1c},\F}=0$, this means $(\omega \otimes \omega)v_{11} = 0$ for all $\omega$. Notice the `polarization' formula $\omega_1 \otimes \omega_2+\omega_2\otimes \omega_1=(\omega_1+\omega_2)\otimes(\omega_1+\omega_2)-\omega_1\otimes\omega_1-\omega_2\otimes\omega_2$, we know that the symmetric 2-tensors of the form $\omega \otimes \omega$ span the space of all symmetric 2-tensors. So from the vanishing condition we conclude that $v_{11}=0$ and consequently $v=0$ in this case. 

Now suppose $\eta_{\mathrm{1c},\F} \neq 0$ and we define $\hat{\eta}_{\mathrm{1c},\F}=\frac{\eta_{\mathrm{1c},\F}}{|\eta_{\mathrm{1c},\F}|}$ then we decompose $\omega$ as $\omega = \epsilon \hat{\eta}_{\mathrm{1c},\F}+(1-\epsilon^2)^{1/2}\omega^\perp$, where $\omega^\perp$ is the unit vector in the direction of the projection of $\omega$ onto $\eta_{\mathrm{1c},\F}^\perp$, the orthogonal complement of $\eta_{\mathrm{1c},\F}$. Then we substitute this into (\ref{kernel_decomp}) to obtain
\begin{align*}
& ( (1+ \frac{|\eta_{\mathrm{1c},\F}|^2}{\xi_{\mathrm{1c},\F}^2+1})^2\epsilon^2 \hat{\eta_{\mathrm{1c},\F}} \otimes \hat{\eta}_{\mathrm{1c},\F})+ (1+ \frac{|\eta_{\mathrm{1c},\F}|^2}{\xi_{\mathrm{1c},\F}^2+1})\epsilon(1-\epsilon^2)^{1/2}
\\&(\hat{\eta}_{\mathrm{1c},\F} \otimes \omega^{\perp} +\omega^\perp \otimes \hat{\eta}_{\mathrm{1c},\F} ) +(1-\epsilon^2)\omega^\perp\otimes\omega^\perp)\cdot v_{11}=0.
\end{align*}
Our condition is that this equation holds for all $\epsilon \in [-1,1]$. Taking $\epsilon = 0$ first yields $(\omega^\perp \otimes \omega^\perp) \cdot v_{11} = 0$. By the same polarization argument as above, cotensors of the form $\omega^\perp \otimes \omega^\perp$ span $\eta_{\mathrm{1c},\F} \otimes \eta_{\mathrm{1c},\F}^\perp$, we conclude that $v_{11}$ is orthogonal to every cotensor in $\eta_{\mathrm{1c},\F}^\perp \otimes \eta_{\mathrm{1c},\F}^\perp$. 
Our second step is to take derivative with respect to $\epsilon$ at $\epsilon=0$ at (\ref{kernel_decomp}), which yields $(\hat{\eta}_{\mathrm{1c},\F} \otimes \omega^\perp + \omega^\perp \otimes \hat{\eta}_{\mathrm{1c},\F}) \cdot v_{11}=0$ for all $\omega^\perp$. Notice that symmetric tensors of the form $(\hat{\eta}_{\mathrm{1c},\F} \otimes \omega^\perp + \omega^\perp \otimes \hat{\eta}_{\mathrm{1c},\F})$ and $\eta_{\mathrm{1c},\F}^\perp \otimes \eta_{\mathrm{1c},\F}^\perp$ together span $(\eta_{\mathrm{1c},\F} \otimes \eta_{\mathrm{1c},\F})^\perp$. Then further taking the second order derivative at $\epsilon=0$ shows that $(\hat{\eta}_{\mathrm{1c},\F} \otimes \hat{\eta}_{\mathrm{1c},\F}) \cdot v_{11}=0$. Combining orthogonal conditions we have we know $v_{11}=0$. Combining (\ref{v_components_relation}), we know that $v=0$, then the non-degeneracy of this principal symbol and thus the ellipticity follows. 
\end{proof}

\subsection{Ellipticity in the combined class}
Next we consider the ellipticity of $N_{h,\F}$ in the combined
operator class $\Psi_{\mathrm{sc},\mathrm{1c},h,\mathcal{F}}$. We
first show the ellipticity of $N_{h,\F}$ as a semiclassical scattering
operator near $\Sigma_{x_0}$ restricted to the kernel of
$\delta_{h,\F}$. Note that this argument is necessary since in \cite{stefanov2018inverting}
only the scattering (not the semiclassical) behavior was considered,
but the addition of the semiclassical behavior does not require any
significant changes.
\begin{prop}
Let $\F>0$ for one forms, and $\F$ is sufficiently large for two tensors and $\Omega_{x_0} = \{ x \leq x_0 \}$ with $x_0$ small.
Then acting on one forms $N_{h,\F}$ is an operator in $\Psi_{\mathrm{sc},h,\mathcal{F}}^{-1,0}$,
while on symmetric 2-tensors, it is an operator in $\Psi_{\mathrm{sc},h,\mathcal{F}}^{-1,2}$. In both cases, it is elliptic on the kernel of $\delta_{h,\F}$ in both scattering algebra and semiclassical sense.
 \label{prop_sc}
\end{prop}
\begin{proof}
The ellipticity of $N_{h,\F}$ as a semiclassical foliation scattering operator on functions is proved in \cite{vasy2020semiclassical}. The process to transplant this result to one forms and 2-tensors is similar to the proof of Proposition \ref{prop_conic}.

As we mentioned in the beginning of this section, the membership of $N_{h,\F}$ near $\Sigma_{x_0}$ in $\Psi_{\mathrm{sc},h,\F}$, expressions of phase and damping factor are given in the proof of Proposition 3.3 of \cite{vasy2020semiclassical} and the argument in Proposition 3.1 of \cite{stefanov2018inverting}.
The remaining part to verify is the power of $x,h$ introduced by our tensorial factor and notice the decay order difference in results. In the 1-form case, as computation below show, we have an extra $h^{-1}(x_0-x)^{-2}$ factor compared with the scalar function case. $h$ is absorbed in the same manner as in our computation near $\Sigma_0$. Thus the membership of $\Psi_{\mathrm{sc},h,\F}^{-1,-2}$ there is changed to be $\Psi_{\mathrm{sc},h,\F}^{-1,0}$. In the 2-tensor case, we have an extra factor $h^{-2}(x_0-x)^{-4}$. An $h^{-1}$ power is absorbed as in the 1-form case, or previous lemmas, and another $h^{-1}$ factor is absorbed by that in the definition of $L$, finally we have membership of $\Psi_{\mathrm{sc},h,\F}^{-1,2}$.

Notice that in the notation \cite{vasy2020semiclassical}, $(x\lambda)$ is interpreted as the component of $\dot{\gamma}$ on $\partial_x$ direction as a whole. 
So the rescaled variable we introduce near $\Sigma_{x_0}$ are
\begin{align}
\tilde{\lambda}:= \frac{x\lambda}{h^{1/2}(x_0-x)}, \, \tilde{t}:= \frac{xt}{h^{1/2}(x_0-x)}.
\label{rescaled_sc}
\end{align}

The metric $g_{\mathrm{sc},\mathrm{1c},h}$ near $\Sigma_{x_0}$, which is a semiclassical scattering metric, applied to (\ref{tangent}) is:
\begin{align*}
g_{\mathrm{sc},h}(x\lambda \partial_x + \omega\partial_y) & = x\lambda \frac{dx}{h^2(x_0-x)^4}+\frac{g_2(\omega\partial_y)}{h(x_0-x)^2}\\
& =  h^{-1/2}(x_0-x)^{-1}(\tilde{\lambda}\frac{dx}{h(x_0-x)^2} + \frac{g_1(\omega\partial_y)}{h^{1/2}(x_0-x)}).
\end{align*}
Similarly, the second factor is
\begin{align*}
\iota_{\dot{\gamma}} = h^{-1/2}(x_0-x)^{-1}( (\tilde{\lambda}+2\alpha\tilde{t})h(x_0-x)^2\partial_x + \omega h^{1/2}(x_0-x)\partial_y).
\end{align*}
Then the counterparts of Lemma \ref{lmn_oneform_infinity} -  \ref{lmn_2tensor_semiclassical} for $N_{h,\F}$ in the semiclassical scattering foliation algebra near $\Sigma_{x_0}$ are proved in the same manner as the proof of Lemma \ref{lmn_oneform_infinity} -  \ref{lmn_2tensor_semiclassical} by replacing $x$ by $(x_0-x)$, $\hat{\lambda}$ by $\tilde{\lambda}$, and $\hat{t}$ by $\tilde{t}$. These results combine to complete our proof. 
\end{proof}

Next we derive the ellipticity of $N_{h,\F}$ with an extra term added in the combined class by combining Proposition \ref{prop_conic} and Proposition \ref{prop_sc}. As before we consider $\Omega_{x_0} = \{ x \leq x_0 \}$ with $x_0$ small.
\begin{prop}
First consider the result about one forms. For $\F>0$, $N_{h,\F}$ is an operator in $\Psi_{\mathrm{sc},\mathrm{1c},h,\mathcal{F}}^{-1,0,-1}(X;_{h,\mathcal{F}}^{\mathrm{sc},\mathrm{1c}}T^*X,_{h,\mathcal{F}}^{\mathrm{sc},\mathrm{1c}}T^*X)$. With suitable choice of $M \in \Psi_{\mathrm{sc},\mathrm{1c},h,\mathcal{F}}^{-3,0,-1}(X)$, the operator 
$$A_{h,\F}:=N_{h,\F}+d_{h,\F}^sM\delta_{h,\F}^s$$
is elliptic in $\Psi_{\mathrm{sc},\mathrm{1c},h,\mathcal{F}}^{-1,0,-1}(X;_{h,\mathcal{F}}^{\mathrm{sc},\mathrm{1c}}T^*X,_{h,\mathcal{F}}^{\mathrm{sc},\mathrm{1c}}T^*X)$ on $\Omega_{x_0}$.

On the other hand, consider symmetric 2-tensors. For $\F$ sufficiently large, $N_{h,\F}$ is an operator in $\Psi_{\mathrm{sc},\mathrm{1c},h,\mathcal{F}}^{-1,2,-1}(X;\mathrm{Sym}_{h,\mathcal{F}}^{2,\mathrm{sc},\mathrm{1c}}T^*X,\mathrm{Sym}_{h,\mathcal{F}}^{2,\mathrm{sc},\mathrm{1c}}T^*X)$.
With suitable choice of $M \in \Psi_{\mathrm{sc},\mathrm{1c},h,\mathcal{F}}^{-3,2,-1}(X;_{h,\mathcal{F}}^{\mathrm{sc},\mathrm{1c}}T^*X,_{h,\mathcal{F}}^{\mathrm{sc},\mathrm{1c}}T^*X)$, the operator 
$$A_{h,\F}:=N_{h,\F}+d_{h,\F}^sM\delta_{h,\F}^s$$
is elliptic in $\Psi_{\mathrm{sc},\mathrm{1c},h,\mathcal{F}}^{-1,2,-1}(X;\mathrm{Sym}_{h,\mathcal{F}}^{2,\mathrm{sc},\mathrm{1c}}T^*X,\mathrm{Sym}_{h,\mathcal{F}}^{2,\mathrm{sc},\mathrm{1c}}T^*X)$.\\
Here, the ellipticity in $\Psi_{\mathrm{sc},\mathrm{1c},h,\mathcal{F}}^{*,*,*}$ includes the ellipticity in the sense of the standard (differential), the scattering (at $\Sigma_{x_0}$) and 1-cusp (at $\Sigma_0$) 
boundary as well as the semiclassical principal symbols. 
In particular, $A_{h,\F}$ is invertible for $h$ sufficiently small.
 \label{prop_combined}
\end{prop}
\begin{proof}
The proofs for one forms and 2-tensors are the same, except for the extra requirement that $\F$ is sufficiently large for 2-tensor cases introduced by the proof of ellipticity in the semiclassical sense.
The membership and ellipticity on the kernel of $\delta_{h,\F}^s$ near two boundaries follows from Proposition \ref{prop_conic} and Proposition \ref{prop_sc}. 
Thus, when we add the term $d_{h,\F}^sM\delta_{h,\F}^s$, where $M$ has positive scalar principal symbol as an operator in $\Psi_{\mathrm{sc},\mathrm{1c},h,\mathcal{F}}^{-3,0,-1}$ for 1-forms and $\Psi_{\mathrm{sc},\mathrm{1c},h,\mathcal{F}}^{-3,2,-1}$ for 2-tensors, then its product with $d_{h,\F}^s,\delta_{h,\F}^s$ is of the same order as $N_{h,\F}$. Since $M$ has positive scalar principal symbol, to show the ellipticity of $A_{h,\F}$, it remains to check the positivity of the principal symbol of $d_{h,\F}^s\delta_{h,\F}^s$.

We only give computation details for the 1-cusp principal symbol near $\Sigma_0$ because the computation near $\Sigma_{x_0}$ is almost the same as that near $\Sigma_0$. 
The membership of $d_{h,\F}^s$ and its adjoint $\delta_{h,\F}^s$ as operator in the combined class is encoded in the discussion in Section \ref{section_ds_1c} and \ref{section_ds_sc}.
Also, according to symbol computation in  Section \ref{section_ds_1c}, on one forms, $d_{h,\F}^s\delta_{h,\F}^s$ has principal symbol
\begin{align*}
\begin{pmatrix}
 \xi_{\mathrm{1c}}-i\F \\  \eta_{\mathrm{1c}} 
 \end{pmatrix} 
\begin{pmatrix}
\xi_{\mathrm{1c}}+i\F && \iota_{\eta_{\mathrm{1c}}}
\end{pmatrix}
=
\begin{pmatrix}
\xi_{\mathrm{1c}}^2+\F^2 && (\xi_{\mathrm{1c}}-i\F)\iota_{\eta_{\mathrm{1c}}}\\
(\xi_{\mathrm{1c}}+i\F)\eta_{\mathrm{1c}} \otimes && \eta_{\mathrm{1c}} \otimes \iota_{\eta_{\mathrm{1c}}}
\end{pmatrix},
\end{align*} 
and on 2-tensors, $d_{h,\F}^s\delta_{h,\F}^s$ has principal symbol
\begin{align}
\scriptsize
\begin{split}
& \begin{pmatrix}
\xi_{\mathrm{1c}}-i\F && 0 \\
\frac{1}{2}\eta_{\mathrm{1c}} \otimes && \frac{1}{2}(\xi_{\mathrm{1c}}-i\F) \\
b_s && \eta_{\mathrm{1c}} \otimes_s
\end{pmatrix}
\begin{pmatrix}
\xi_{\mathrm{1c}}+i\F && \frac{1}{2}\iota_{\eta_{\mathrm{1c}}} && \la b_s, \cdot \ra \\
0 && \frac{1}{2}(\xi_{\mathrm{1c}}-i\F) && \iota_{\eta_{\mathrm{1c}}}
\end{pmatrix}
=\\
&  \begin{pmatrix}
\xi_{\mathrm{1c}}^2+\F^2 &&  \frac{1}{2}(\xi_{\mathrm{1c}}-i\F)\iota_{\eta_{\mathrm{1c}}} && (\xi_{\mathrm{1c}}-i\F)\la b_s, \cdot \ra \\
\frac{1}{2}(\xi_{\mathrm{1c}}+i\F)\eta_{\mathrm{1c}} \otimes && \frac{1}{4}((\eta_{\mathrm{1c}} \otimes \cdot)\iota_{\eta_{\mathrm{1c}}} + (\xi_{\mathrm{1c}}^2+\F^2) ) &&  \frac{1}{2} \eta_{\mathrm{1c}} \otimes \la b_s, \cdot \ra+\frac{1}{2}(\xi_{\mathrm{1c}}-i\F)\iota_{\eta_{\mathrm{1c}}}^s\\
(\xi_{\mathrm{1c}}-i\F)b_s  && \frac{1}{2}b_s\iota_{\eta_{\mathrm{1c}}} + \frac{1}{2}(\xi_{\mathrm{1c}}+i\F)\eta_{\mathrm{1c}}\otimes_s &&
b_s\la b_s,\cdot \ra + \eta_{\mathrm{1c}}\otimes_s \iota_{\eta_{\mathrm{1c}}}
\end{pmatrix},
\end{split}
\label{ddelta_symbol}
\end{align}
where the last matrix in fact is a $4\times4$ matrix of blocks, where the thrid column is the same as the second column and omitted. When $b_s=0$, standard linear algebra gives lower bound independent of $\F$ in terms of $|(\xi_{\mathrm{1c}},\eta_{\mathrm{1c}})|$. Thus, rescaling $\xi_{\mathrm{1c}},\eta_{\mathrm{1c}}$ by $\F^{-1}$ as in the proof of Lemma \ref{lmn_2tensor_semiclassical}, we can absorb terms involving $b_s,b_s^*$, since their total power of $\F,\xi_{\mathrm{1c}},\eta_{\mathrm{1c}}$ is lower.

To summarize, this extra term $d_{h,\F}^sM\delta_{h,\F}^s$, with $M$ suitably large and positive gives the ellipticity of $A_{h,\F}$, without changing its action on the kernel of $\delta_{h,\F}^s$.
\end{proof}

\subsection{Proof of Theorem \ref{thm_conic}}
\begin{proof}
Proposition \ref{prop_combined} implies the invertibility of
$A_{h,\F}$, and thus, under the gauge condition, the injectivity of
$N_{h,\F}$, which in turn implies the injectivity of $Ie^{\frac{\F\Phi}{h}}$ restricted to tensors in the kernel of $\delta_{h,\F}^s$ with sufficient decay and hence the result of Theorem \ref{thm_conic}.
\end{proof}

\begin{rmk}
  Following up on Remark~\ref{rmk_general_p_0}, see also analogous
  arguments in \cite{zachos2022inverting}, when we replace $x$ by $x^p$, then the rescaled variables should be taken as $\hat{\lambda}=\lambda/(h^{1/2}x^p),\hat{t}=t/(h^{1/2}x^p)$, and similarly the weight for conjugation is $\Phi=-\frac{1}{2px^{2p}}$ near $\Sigma_0$. The 1-cusp algebra is changed accordingly: the construction is completely the same, but with the defining function $x$ and the smooth structure replaced by that of $x^p$.  


Now write $\xi_{\mathrm{1c},p},\eta_{\mathrm{1c},p}$ for the dual variables of the 1-cusp cotangent bundle constructed using $x^p$ as the boundary defining function. Covectors are
\begin{align*}
\xi_{\mathrm{1c},p}\frac{dx^p}{x^{3p}} + \eta_{\mathrm{1c},p} \frac{dy}{x^p}  = p \xi_{\mathrm{1c},p}\frac{dx}{x^{2p+1}}+\eta_{\mathrm{1c},p}\frac{dy}{x^p}. 
\end{align*}
The phase (\ref{phase_conic_1}) becomes (using our new $\hat{t},\hat{\lambda}$)
\begin{align*}
& p\xi_{\mathrm{1c},p}( \hat{\lambda}\hat{t}+\alpha(x,y,h^{1/2}x^p\hat{\lambda},\omega)\hat{t}^2+h^{1/2}x^p\hat{t}^3\Gamma^{(1)}(x,y,h^{1/2}x^p\hat{\lambda},\omega,h^{1/2}x^p\hat{t})) \
\\& + \eta_{\mathrm{1c},p} \cdot (\omega \hat{t} + h^{1/2}x^p\hat{t}^2 \Gamma^{(2)}(x,y,h^{1/2}x^p\hat{\lambda},\omega,h^{1/2}x^p\hat{t}) ).
\end{align*}
Similarly, the damping factor (\ref{damp}) becomes
\begin{align*}
\hat{\lambda}\hat{t} + \alpha(x,y,h^{1/2}x^p\hat{\lambda},\omega)\hat{t}^2+h^{1/2}x^p\hat{t}^3\Gamma^{(1)}(x,y,h^{1/2}x^p\hat{\lambda},\omega,h^{1/2}x^p\hat{t}).
\end{align*}
Then the arguments afterwards for symbol computation and ellipticity in various senses go through as before after repalcing quantities as mentioned here, and the $p$ coefficient produced does not affect the argument.
\label{rmk_general_p}
\end{rmk}

\section{The gauge condition}
\label{sec_gauge}
\subsection{The gauge condition, conic 2-tensor}
We still need to show that we can arrange the gauge condition 
\begin{align*}
\delta_{h,\F}^sf_{h,\F}=0,
\end{align*}
where
\begin{align*}
f_{h,\F} = e^{-\frac{\F\Phi}{h}} f.
\end{align*} 
Recall that the freedom we have is adding to $f$ a term of the form
$d^sv$ with $v$ decaying suffiicently fast at $\partial \overline{M}$. This is equivalent to adding to $f_{h,\F}$ a term of the form $d_{h,\F}^s v_{h,\F}$ with $v_{h,\F} = e^{-\frac{\F\Phi}{h}} v$, where $d_{h,\F}^s = e^{-\frac{\F\Phi}{h}} d^s  e^{\frac{\F\Phi}{h}}$ is the adjoint of $\delta_{h,\F}$, and the modified Laplacian is their product:
\begin{align*}
\Delta_{h,\F,s} = \delta_{h,\F}^s d_{h,\F}^s.
\end{align*}
Notice that if $v_{h,\F}$ is in any 1-cusp Sobolev space (of
sufficient regularity, if one wants pointwise statements) near $\partial\overline{M}$ then $v$ is
actually Gaussian decaying, so $d^sv$ is indeed in the kernel of the
X-ray transform. Note also that there is no decay needed at the
artificial boundary: the geodesics on which our modified normal
operator puts a positive weight do not intersect it! However, we must
of course make sure that the added potential term leaves us is the
correct function space.

The modified solenoidal ($\mathcal{S}$) and potential ($\mathcal{P}$) projections acting on a function or one form $\phi$ are given by
\begin{align*}
& \mathcal{S}_{h,\F}\phi = \phi - d_{h,\F}^s\Delta_{h,\F,s}^{-1}\delta_{h,\F}^s\phi,\\
& \mathcal{P}_{h,\F}\phi = d_{h,\F}^sQ_{h,\F}\phi, \, Q_{h,\F}\phi = \Delta_{h,\F,s}^{-1}\delta_{h,\F}^s\phi.
\end{align*} 
$Q_{h,\F}\phi$ is vanishing at $\partial_{\mathrm{int}}{\Omega}$ because of the boundary condition for $\Delta_{h,\F,s}$ and its solution operator. Thus $\mathcal{P}_{h,\F}\phi$ is in the range of $d_{h,\F}^s$ applied to functions or one forms vanishing at $\partial_{\mathrm{int}}{M}$. Further, $\mathcal{S}_{h,\F}\phi$ is in the kernel of $\delta_{h,\F}$.
\begin{align*}
\delta_{h,\F}\mathcal{S}_{h,\F}\phi = \delta_{h,\F}^s\phi- \delta_{h,\F}^s d_{h,\F}^s\Delta_{h,\F,s}^{-1}\delta_F^s\phi=0.
\end{align*} 
Thus the remaining task is to justify the definition of $Q_{h,\F}$ by
checking the invertibility of $\Delta_{h,\F,s}$ on the mixed 1-cusp
(at $\partial\overline{M}$)/scattering (at the artificial boundary) Sobolev spaces.

\subsection{The invertibility of $\Delta_{h,\F,s}$}
In this section we prove the ellipticity of the modified Laplacian $\Delta_{h,\F,s}$, from which the invertibility follows by the taking semiclassical limit. This, in combination with previous discussion, proves Theorem \ref{thm_gauge}.
\begin{lmn}
For $\F>0$, we have $\Delta_{h,\F,s}=\delta_{h,\F}^sd_{h,\F}^s \in
\Psi_{\mathrm{sc},\mathrm{1c},h,\mathcal{F}}^{2,0,0}(X)$, and for
$\F>0$ sufficiently large it is jointly elliptic in the sense of the
standard (differential), the 1-cusp (at $\Sigma_0$) and scattering (at $\Sigma_{x_0}$) boundary as well as the semiclassical principal symbols. In particular, it is invertible for $h$ sufficiently small.
\label{lmn_Delta_invert}
\end{lmn}
\begin{proof}
Membership of $\Delta_{h,\F,s}$ near both boundaries follows from the operator properties of the $d_{h,\F}^s$ discussed in Section \ref{section_ds_1c} and \ref{section_ds_sc}  and consequently its adjoint $\delta_{h,\F}^s$. Near the scattering boundary $\Sigma_{x_0}$, it is in $\mathrm{Diff}_{h,\mathrm{sc}}^{2,0}(X)$; while near the 1-cusp boundary $\Sigma_0$, it is in $\mathrm{Diff}_{\mathrm{1c}}^{2,0}(X)$. The bundle valued version can be derived in the same manner. Combining these facts, we know that in our new operator class, $\Delta_{h,\F,s}$ lies in $\Psi_{\mathrm{sc},\mathrm{1c},h,\mathcal{F}}^{2,0,0}$.

Next we consider the ellipticity near $\Sigma_0$ as an 1-cusp operator. Recall Proposition \ref{prop_delta_symbol}, the principal symbol of $\Delta_{h,\F,s}$ is
\begin{align*}
\begin{pmatrix}
\xi_{\mathrm{1c}}^2+\F^2+\frac{1}{2}\eta_{\mathrm{1c}}^2 && \frac{1}{2}(\xi_{\mathrm{1c}}-i\F)\iota_{\eta_{\mathrm{1c}}} \\
\frac{1}{2}(\xi_{\mathrm{1c}}+i\F)\eta_{\mathrm{1c}} \otimes && \frac{1}{2}(\xi_{\mathrm{1c}}^2+\F^2)+ \iota_{\eta_{\mathrm{1c}}}^s \eta_{\mathrm{1c}}\otimes_s
\end{pmatrix}
 + \begin{pmatrix}
\la b_s,\cdot \ra b_s && \la b_s,\cdot \ra \eta_{\mathrm{1c}}\otimes_s \\
\iota_{\eta_{\mathrm{1c}}}^sb_s && 0
\end{pmatrix},
\end{align*}
where the term $\frac{1}{2}(\xi_{\mathrm{1c}}^2+\F^2)+\iota_{\eta_{\mathrm{1c}}}^s \eta_{\mathrm{1c}}\otimes_s$ represents a block with components $\frac{1}{2}(\xi_{\mathrm{1c}}^2+\F^2)+\frac{1}{2}|\eta_{\mathrm{1c}}|^2\delta_{ij}+\frac{1}{2}\eta_{\mathrm{1c},i}\eta_{\mathrm{1c},j}$.
We apply the same argument in the proof of Lemma \ref{lmn_2tensor_semiclassical} to absorb the second term by taking $\F>0$ sufficiently large since all those terms have lower order in terms of the total power of $\xi_{\mathrm{1c}},\eta_{\mathrm{1c}}$ and $\F$. This proves the ellipticity near $\Sigma_0$ in the 1-cusp and semiclassical sense.

The ellipticity of $\Delta_{h,\F,s}$ in the scattering and semiclassical sense near $\Sigma_{x_0}$ follows from the same argument, but using Proposition \ref{prop_delta_symbol_sc} instead of Proposition \ref{prop_delta_symbol}. The invertibility follows the standard way by constructing parametrix and take semiclassical limit.
\end{proof}

In particular, this proves Theorem~\ref{thm:gauging}.

\begin{proof}[Proof of Theorem~\ref{thm:gauging}]
This is an immediate consequence of the discussion before Lemma \ref{lmn_Delta_invert} and Lemma \ref{lmn_Delta_invert} itself.
\end{proof}

\bibliographystyle{plain}
\bibliography{bib_inverse}

\begin{thebibliography}{10}

\bibitem{Guillarmou-Lassas-Tzou:Conic}
Colin Guillarmou, Matti Lassas, and Leo Tzou.
\newblock X-ray transform in asymptotically conic spaces.
\newblock {\em Int. Math. Res. Not. IMRN}, (5):3918--3976, 2022.

\bibitem{Guillarmou-Mazzucchelli-Tzou:Conjugate}
Colin Guillarmou, Marco Mazzuchelli, and Leo Tzou.
\newblock Asymptotically {E}uclidean metrics without conjugate points are flat.
\newblock {\em Preprint, arXiv:1909.01488}, 2019.

\bibitem{RBMZw}
R.~B. Melrose and M.~Zworski.
\newblock Scattering metrics and geodesic flow at infinity.
\newblock {\em Inventiones Mathematicae}, 124:389--436, 1996.

\bibitem{Melrose1994}
Richard~B. Melrose.
\newblock Spectral and scattering theory for the {L}aplacian on asymptotically
  {E}uclidian spaces.
\newblock In {\em Spectral and scattering theory ({S}anda, 1992)}, volume 161
  of {\em Lecture Notes in Pure and Appl. Math.}, pages 85--130. Dekker, New
  York, 1994.

\bibitem{paternain2019geodesic}
Gabriel~P. Paternain, Mikko Salo, Gunther Uhlmann, and Hanming Zhou.
\newblock The geodesic {X}-ray transform with matrix weights.
\newblock {\em Amer. J. Math.}, 141(6):1707--1750, 2019.

\bibitem{stefanov2018inverting}
Plamen Stefanov, Gunther Uhlmann, and Andr\'{a}s Vasy.
\newblock Inverting the local geodesic {X}-ray transform on tensors.
\newblock {\em J. Anal. Math.}, 136(1):151--208, 2018.

\bibitem{uhlmann2016inverse}
Gunther Uhlmann and Andr{\'a}s Vasy.
\newblock The inverse problem for the local geodesic ray transform.
\newblock {\em Invent. Math.}, 205(1):83--120, 2016.

\bibitem{uhlmann2016journey}
Gunther Uhlmann and Hanming Zhou.
\newblock Journey to the center of the earth.
\newblock {\em arXiv preprint arXiv:1604.00630}, 2016.

\bibitem{vasy2020semiclassical}
Andr{\'a}s Vasy.
\newblock A semiclassical approach to geometric {X}-ray transforms in the
  presence of convexity.
\newblock {\em arXiv preprint arXiv:2012.14307}, 2020.

\bibitem{zachos2022inverting}
Andr{\'a}s Vasy and Evangelie Zachos.
\newblock The {X}-ray transform on asymptotically conic spaces, 2022.

\bibitem{Zachos:Thesis}
Evangelie Zachos.
\newblock {\em The {X}-ray transform on asymptotically {E}uclidean spaces}.
\newblock PhD thesis, Stanford University, 2020.

\end{thebibliography}

\end{document}